\documentclass[12pt]{article}
\usepackage{latexsym}
\usepackage{amsmath,amsfonts}
\usepackage{amssymb}
\usepackage{epsfig}
\usepackage{color}
\usepackage[subnum]{cases}
\usepackage{appendix}
\usepackage{stmaryrd}
\usepackage{cancel}

\usepackage{mathrsfs}
\newenvironment{proof}[1][Proof]{\noindent\textbf{#1.} }{\ \rule{0.5em}{0.5em}}

\newcommand{{\Omegapm}}{\Omega^{i,e}}
\newcommand{{\Oi}}{\Omega_i}
\newcommand{{\Oo}}{\Omega_0}

\newcommand{{\Gammai}}{\Gamma_i}
\newcommand{{\Gammao}}{\Gamma_0}
\newcommand{{\gDpm}}{\gamma_{D;i}^{\pm}}
\newcommand{{\gNpm}}{\gamma_{N;i}^{\pm}}
\newcommand{{\gDp}}{\gamma_{D;i}^{+}}
\newcommand{{\gNp}}{\gamma_{N;i}^{+}}
\newcommand{{\gDm}}{\gamma_{D;i}^{-}}
\newcommand{{\gNm}}{\gamma_{N;i}^{-}}
\newcommand{{\Oipm}}{\Omega_i^{\pm}}
\newcommand{{\Oip}}{\Omega_i^{+}}
\newcommand{{\Oim}}{\Omega_i^{-}}

\newcommand{{\ui}}{u_i}
\newcommand{{\uo}}{u_0}
\newcommand{{\gi}}{\gamma_i}
\newcommand{{\gDop}}{\gamma_{D;0}^{+}}
\newcommand{{\gNop}}{\gamma_{N;0}^{+}}

\newcommand{\kj}{\kappa_j}
\newcommand{\uDi}{u_{D;i}}
\newcommand{\uNi}{u_{N;i}}
\newcommand{\uNj}{u_{N;j}}

\newcommand{\vDj}{v_{D;j}}

\newcommand{\vNj}{v_{N;j}}

\newcommand{\uDj}{u_{D;j}}
\newcommand{\uDl}{u_{D;\ell}}
\newcommand{\uNl}{u_{N;\ell}}

\newcommand{\uDo}{u_{D;0}}
\newcommand{\uNo}{u_{N;0}}

\newcommand{\uDoi}{u^i_{D;0}}
\newcommand{\uNoi}{u^i_{N;0}}

\newcommand{\uDoj}{u^j_{D;0}}
\newcommand{\uNoj}{u^j_{N;0}}
\newcommand{\vDoj}{v^j_{D;0}}
\newcommand{\vNoj}{v^j_{N;0}}
\newcommand{\uDol}{u^\ell_{D;0}}
\newcommand{\uNol}{u^\ell_{N;0}}

\newcommand{\uDoo}{u^0_{D;0}}
\newcommand{\uNoo}{u^0_{N;0}}
\newcommand{\vDoo}{v^0_{D;0}}
\newcommand{\vNoo}{v^0_{N;0}}

\newcommand{\Vj}{V_j}
\newcommand{\Kj}{{K}_j}
\newcommand{\Kpj}{{K}^\prime_j}
\newcommand{\Wj}{{W}_j}
\newcommand{\Vo}{{V}_0}
\newcommand{\Ko}{{K}_0}

\newcommand{\Wo}{{W}_0}

\newcommand{\Vojl}{{V}_0^{j,\ell}}
\newcommand{\Kojl}{{K}_0^{j,\ell}}
\newcommand{\Kpojl}{\left(K^\prime\right)_0^{j,\ell}}
\newcommand{\Wojl}{{W}_0^{j,\ell}}
\newcommand{\Vojo}{{V}_0^{j,0}}
\newcommand{\Kojo}{{K}_0^{j,0}}
\newcommand{\Kpojo}{\left(K^\prime\right)_0^{j,0}}
\newcommand{\Wojo}{{W}_0^{j,0}}

\newcommand{\Vool}{{V}_0^{0,\ell}}
\newcommand{\Kool}{{K}_0^{0,\ell}}
\newcommand{\Kpool}{\left(K^\prime\right)_0^{0,\ell}}
\newcommand{\Wool}{{W}_0^{0,\ell}}
\newcommand{\Vooo}{{V}_0^{0,0}}
\newcommand{\Kooo}{{K}_0^{0,0}}
\newcommand{\Kpooo}{\left(K^\prime\right)_0^{0,0}}
\newcommand{\Wooo}{{W}_0^{0,0}}

\renewcommand{\Xi}{X_i}

\newcommand{\uvj}{\underline{v}_j}

\newcommand{\uuo}{\underline{u}_0}
\newcommand{\uvo}{\underline{v}_0}
\newcommand{\ugo}{\underline{g}_0}

\newcommand{\uuuno}{\underline{u}_1}
\newcommand{\uuq}{\underline{u}_q}
\newcommand{\uvuno}{\underline{v}_1}
\newcommand{\uvq}{\underline{v}_q}
\def\doubleunderline#1{\underline{\underline{#1}}}
\newcommand{\uuu}{\doubleunderline{u}}
\newcommand{\uuC}{\doubleunderline{C}}

\newcommand{\uuv}{\doubleunderline{v}}

\newcommand{\Om}{\Omega}

\newcommand{\no}{{\bf n}_0}
\renewcommand{\ni}{{\bf n}_i}
\newcommand{\R}{\mathbb R}
\newcommand{\C}{\mathbb C}

\newcommand{\x}{{\mathbf x}}
\newcommand{\y}{{\mathbf y}}

\newcommand{\bg}{{\mathbf{g}}}

\newtheorem{lemma}{Lemma}[section]
\newtheorem{theorem}[lemma]{Theorem}

\newtheorem{proposition}[lemma]{Proposition}

\newtheorem{remark}[lemma]{Remark}

\begin{document}

\title{A spectral boundary element method for acoustic interference problems 
\thanks{We acknowledge that the present research has been supported by GNCS-INDAM 2024 research program \emph{``Tecniche nu\-me\-ri\-che efficienti per problemi differenziali avanzati:
	BEM e paradigma isogeometrico"}.}}
\author{S. Falletta\footnote{Dipartimento di Scienze Matematiche, Politecnico di Torino, Italy. Email: silvia.falletta@polito.it}, S. Sauter\footnote{Institut für Mathematik, Universität Zürich, Zürich. Email: stas@math.uzh.ch},
}
\date{}
\maketitle

\abstract{}
In this paper we consider high-frequency acoustic transmission problems with
jumping coefficients modelled by Helmholtz equations. The solution then is
highly oscillatory and, in addition, may be localized in a very small vicinity
of interfaces (whispering gallery modes). For the reliable numerical
approximation a) the PDE is tranformed in a classical single trace integral
equation on the interfaces and b) a spectral Galerkin boundary element method
is employed for its solution. We show that the resulting integral equation is
well posed and analyze the convergence of the boundary element method for the
particular case of concentric circular interfaces. We prove a condition on the
number of degrees of freedom for quasi-optimal convergence. Numerical
experiments confirm the efficiency of our method and the sharpness of the theoretical estimates.


\section{Introduction}
The numerical simulation of high-frequency scattering problems modelled by
Helmholtz equations is challenging and a topic of active research in numerical
analysis and scientific computing. A large wavenumber has the impact that the
arising sesquilinear form in a variational formulation is highly indefinite
causing severe stability issues for the numerical discretization (see, e.g.,
\cite{Schatz74}, \cite{Ihlenburg}, \cite{BabuskaSauter}, \cite{Sauter2005},
\cite{zhu-wu12b}). If the coefficients in the Helmholtz equation (refractive
index, wave speed, material density) are non-constant additional physical
phenomena may arise such as the localization of waves, interference etc. which
makes its numerical solution notoriously hard. In this paper we consider the
Helmholtz problem with coefficients being constant on isolated subdomains with
smooth boundary. It is known that there are geometric configurations of the
subdomains where, for certain critical (quasi-resonant) relations between the
coefficients and the global wavenumber, the wave strongly \textit{localizes}
at the interfaces between different subdomains. This phenomenon, also known as
whispering gallery modes (see, e.g., \cite{ilchenko2006optical},
\cite{balac2021asymptotics}), is well studied in the literature for certain
model cases such as spherical symmetric domains and coefficients (see
\cite{MoiolaSpence_resonance_2017}, \cite{Torres_Sauter_2}). 

In this paper we propose the following numerical solution method.

1) The fact that the wave may localize around an interface and then decays
exponentially away from the interface makes the use of a skeleton integral
equation appealing: the partial differential equations is transformed as a
non-local integral equation on the domain skeleton, i.e., the union of the
interfaces and the domain boundary. This leads to a reduction of the dimension
of the problem by $1$ since the integral equation has to be solved on the
lower-dimensional interface and, in turn, the number of numerical degrees of
freedom is reduced substantially. There are many ways of transforming boundary
value problems to integral equations and we have chosen the classical direct
single trace formulation of the first kind as in \cite{vonPetersdorff89},
\cite{Hiptmair_multiple_trace}, \cite{florian2023skeleton}, \cite{GrHiSa_pw_Lip} so that
the stability properties of the original transmission PDE can be inherited to
the integral equation. In this way, the problem of spurious resonances arising
in indirect formulations (see, e.g., \cite[Sec. 7.7]{ChenZhou}) can be avoided
as well as the use of combined field integral equations (cf. see, e.g.,
\cite{hiptmair2003coercive}).

2) We assume that the interface is analytic and propose a spectral Galerkin
discretization for the discretization of the skeleton integral equations. This
is particularly appealing for circular interfaces, where the eigensystem of the
integral operators on these interfaces are known; we refer, e.g., to
\cite{Amini}, \cite{SaSchw1}, \cite{jerez2020high}, \cite{jerez2022spectral},
\cite{chandler2007galerkin}, \cite{ecevit2022spectral} for various aspects and
applications of the spectral Galerkin boundary element method.
% new
Spectral Galerkin methods can also be applied to non-circular interfaces,
e.g., by using mapped spherical harmonics or the $p$-version of the boundary
element method. For a concentric circular interface, we prove that optimal
convergence starts if the number of freedoms is proportionally to the
wavenumber. As a nice surprise, this condition for our Helmholtz problem with
\textit{piecewise constant coefficients} and possible interference is
qualitatively the same as the condition on the order of the (spectral) method
for the Helmholtz equation with\textit{ constant coefficients}, see, e.g.,
\cite[Thm. 3.17 (setting $h=O\left(  1\right)  $ and $p\gtrsim k$%
)]{MelenkLoehndorf} for boundary integral equations and \cite[Thm.
3.9]{hiptmair-moiola-perugia09b}, \cite[Thm. 5.5 (setting $h=O\left(
1\right)  $ and $p\gtrsim k$)]{MelenkSauterMathComp} for PDE formulations of
the Helmholtz problem. We emphasize that a stability and convergence analysis
of our method for \textit{non-circular} interfaces is still an open problem
for future research.\bigskip

The combination of the above mentioned two techniques results in a discrete linear system
whose dimension is much smaller compared to the direct discretization of the
PDE or to the $h$-version of the boundary element method (see, e.g.,
\cite{SauterSchwab2010}) and could allow for the reliable simulation of
high-frequency Helmholtz problems even for localized waves. 

We analyze theoretically the convergence of the spectral Galerkin method for
the special case of a transmission problem on a disc with concentric
interface. The main result can be phrased as follows: if the ratio of the
(spectral) degrees of freedom and the wavenumber is smaller than a fixed
constant (depending linearly on the size of material parameters) the
convergence is quasi-optimal even in quasi-resonant cases. We have performed
numerical experiments which confirm the efficiency of our method and the
sharpness of the theoretical estimates.

The paper is structured as follows. In Section \ref{sec:standard} the acoustic
transmission model problem is formulated along the definition of the involved
operators and the assumptions on parameters. Section \ref{sec:BIE} is devoted to
the derivation of the well-posed single trace formulation of this problem on
the interfaces and domain boundary. In Section \ref{sec:concentric_cirles} we consider the
specific case of two concentric circular interfaces and introduce the spectral
Galerkin boundary element method. The chosen boundary element basis
diagonalizes all involved operators and explicit formulas for its diagonal
entries are given. This particular setting allows us to derive regularity
estimates for the traces and jumps of the exact solution depending on the
smoothness of the given boundary data. In Section \ref{sec:numerical_results} we report on the
results of numerical experiments which illustrate the findings in Section
\ref{sec:concentric_cirles}. Finally, some estimates of Bessel and exponential functions are
presented in Appendix \ref{sec:appendix}.

\section{The model problem}\label{sec:standard}
\begin{figure}[h!]
	\centering
	\includegraphics[width=0.3\textwidth]{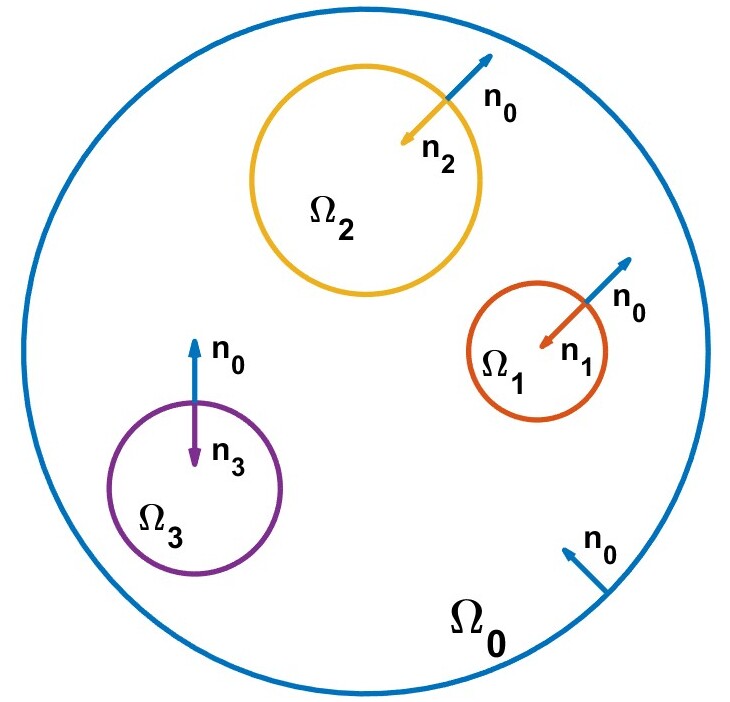}
	\caption{Model problem geometrical setting.}
	\label{fig:pb_domain}
\end{figure}

Let $\Oi\subset \R^2$, with $i = 1,\cdots, q$, be open, bounded rigid and disjoint domains, whose boundaries $\Gammai = \partial\Oi$ are assumed to be closed and smooth curves. Let $\Omega^e = \mathbb{R}^2\setminus (\cup_{i = 1}^q \overline{\Oi})$ and denote by $\Omega_0$ the bounded domain obtained by delimiting the unbounded region $\Omega^e$ by a close and regular  contour $\Gamma_0$. It then results that $\partial\Omega_0 = \Gammao\cup\left(\cup_{i=1}^q \Gammai\right)$.
%Then, we set $\Oip = \Oi$ and $\Oim := \R^2 \setminus \overline{\Oi}$. 
Let $\ni$ be the inner normal vector field for the domain $\Oi$, $i = 0,\cdots,q$ (see Figure \eqref{fig:pb_domain}).
%%%%

{\color{black}Given a real number $s\geq 0$, and a general domain $\mathcal{O}\subset \mathbb{R}^2$, 
we denote by $\| \cdot \|_{H^s(\mathcal{O})}$ 
%and $| \cdot |_{s,\mathcal{O}}$ respectively, 
the norm 
%and seminorm 
of the usual
Sobolev space $H^s(\mathcal{O})$ (cf. \cite{Adams}). Also, 
%we use the convention $L^2(\mathcal{O}) := H^0(\mathcal{O})$, and 
for all $t \in (0,1]$ we let $H^{-t}(\partial\mathcal{O})$ be the dual of $H^{t}(\partial\mathcal{O})$ with respect to the pivot space $L^2(\partial\mathcal{O})$.}
%%%%%

Let $\kappa \in \C\setminus\{0\}$, with {\color{black}$\mathfrak{\operatorname{Im}}(\kappa) \geq 0$}, denote the wavenumber and let $n_i >0$ be real numbers.
We consider the following Helmholtz transmission problem: find $u(\x)$, $\x\in\Om = \cup_{i=0}^q \Oi$ such that, denoting by $\ui = u_{|_{\Oi}}$,

\begin{subnumcases}{}
	\Delta \ui(\x)+\kappa^{2}n_i \ui(\x)= 0 & $\x\in\Oi$, \label{transmission_problem_1}
	\\
	\gamma_{D;i}u_i(\x) = \gamma_{D;0}^i u_0(\x), & $\x\in \Gammai$, $i = 1,\cdots,q$, \label{transmission_problem_2}
	\\ 
	\gamma_{N;i}u_i(\x) = -\gamma_{N;0}^i u_0(\x), & $\x\in \Gammai$, $i = 1,\cdots,q$, \label{transmission_problem_3}
	\\
	(\gamma_{N;0}^0 u_0)(\x) - \text{DtN}(\gamma_{D;0}^0 \uo)(\x) = g(\x) & $\x\in \Gammao$, \label{transmission_problem_4}
	%\frac{\partial u_0(\x)}{\partial\no} - \text{DtN}(\gDop\uo) = g(\x) & $\x\in \Gammao$ \label{transmission_problem_3}
\end{subnumcases}
where $\gamma_{D;i}$ and $\gamma_{D;0}^i$ in \eqref{transmission_problem_2} are the {\color{black} interior} Dirichlet trace operators
\begin{equation}\label{eq:trace_D}
\begin{aligned}
	&\gamma_{D;i} : H^1(\Oi) \longrightarrow H^{1/2}(\Gammai), & & \qquad i = 1,\cdots,q\\
	&\gamma_{D;0}^i : H^1(\Omega_0) \longrightarrow H^{1/2}(\Gammai), & & \qquad i = 0,\cdots,q,
\end{aligned}
\end{equation}
and $\gamma_{N;i}$ and $\gamma_{N;0}^i$ in \eqref{transmission_problem_3} are the {\color{black} interior} Neumann trace operators
\begin{equation}\label{eq:trace_N}
\begin{aligned}
	&{\color{black}\gamma_{N;i}} : H^1(\Oi;\Delta) \longrightarrow H^{-1/2}(\Gammai), & & \qquad i = 1,\cdots,q\\
	&\gamma_{N;0}^i : H^1(\Omega_0;\Delta) \longrightarrow H^{-1/2}(\Gammai), & & \qquad i = 0,\cdots,q,
	\end{aligned}
\end{equation}
	{\color{black} $H^1(\mathcal{O};\Delta)$ being the subspace of $H^1(\mathcal{O})$ made of functions with square-integrable (distributional) Laplacian: $H^1(\mathcal{O};\Delta) = \{v\in H^1(\mathcal{O}) \, : \, \Delta v \in L^2(\mathcal{O})\}$}.
Finally, in the boundary condition \eqref{transmission_problem_4} for the accurate representation of the unbounded domain, the $\text{DtN}$ denotes the Dirichlet-to-Neumann operator that maps the Dirichlet trace of $u_0$ on $\Gammao$ to its Neumann trace. {\color{black}The well-posedness of problem \eqref{transmission_problem_1}--\eqref{transmission_problem_4} is proved, e.g., in \cite[Sec.
2]{vonPetersdorff89}, \cite[Sec. 2.2]{GrHiSa_pw_Lip}. }

\

In what follows, to simplify the notation, for $u_i$ sufficiently regular, we will denote by
\begin{equation}\label{eq:trace_notation}
\begin{aligned}
&\uDi = \gamma_{D;i} u_i = u_i{_{|_{\Gammai}}}, &  & \qquad i = 1,\cdots,q, \\
&\uDoi = \gamma_{D;0}^i u_0 = u_0{_{|_{\Gammai}}}, & &\qquad i = 0,\cdots,q, \\
&\uNi = \gamma_{N;i} u_i = \partial_{\ni} \ui,  & &\qquad i = 1,\cdots,q,\\
&\uNoi = \gamma_{N;0}^i u_0 = \partial_{\no(\x)} \uo, & &\qquad  \x\in\Gammai,  \, i = 0,\cdots,q. \\
\end{aligned}
\end{equation}

\section{The bounday integral equation formulation}\label{sec:BIE}

Let $G_{\kj}(\x,\y) = \frac{\operatorname*{i}}{4} H^{(1)}_0(\kj\|\x-\y\|)$, with $\kj = \kappa\sqrt{n_j}$, $j = 0,\cdots,q$,  be the fundamental solution of the Helmholtz problem, $H^{(1)}_0$ being the Hankel function of first kind and order zero.
We define the single and double layer potentials related to the domain $\Omega_j$ by
\[%
\begin{array}
	[c]{cl}%
	S_{\ell}^{j}:H^{-1/2}\left(  \Gamma_{j}\right)  \rightarrow
	H_{\operatorname*{loc}}^{1}\left(  \mathbb{R}^{2}\right)   & \left(
	S_{{\ell}}^{j}\varphi\right)  \left(  \mathbf{x}\right)  :=\int
	_{\Gamma_{j}}G_{\kappa_{\ell}}\left(  \mathbf{x},\mathbf{y}\right)
	\varphi\left(  \mathbf{y}\right)  \text{d}\mathbf{y},\\
	D_{\ell}^{j}:H^{1/2}\left(  \Gamma_{j}\right)  \rightarrow
	H_{\operatorname*{loc}}^{1}\left(  \mathbb{R}^{2}\setminus\Gamma_j\right)   & \left(
	D_{{\ell}}^{j}\varphi\right)  \left(  \mathbf{x}\right)  :=\int
	_{\Gamma_{j}}\left\langle  \mathbf{n}_{j}\left(  \mathbf{y}\right)
	,\nabla_{\mathbf{y}}G_{\kappa_{\ell}}\left(  \mathbf{x},\mathbf{y}\right)
	\right\rangle \varphi\left(  \mathbf{y}\right)  \text{d}\mathbf{y}.
\end{array}
\]
This gives rise to the following boundary integral operators: for sufficiently
smooth $\varphi$ and $\psi$ we define

\begin{itemize}
	\item Single layer operators:%
	\begin{align*}
		V_{0}^{j,\ell}\varphi\left(  \mathbf{x}\right)   &  :=\int_{\Gamma_{\ell}%
		}G_{\kappa_{0}}\left(  \mathbf{x},\mathbf{y}\right)  \varphi\left(
		\mathbf{y}\right)  \text{d}\mathbf{y\quad x}\in\Gamma_{j}\quad j,\ell \in\left\{
		0,\cdots,q\right\}  ,\\
		V_{j}\varphi\left(  \mathbf{x}\right)   &  :=\int_{\Gamma_{j}}G_{\kappa_{j}%
		}\left(  \mathbf{x},\mathbf{y}\right)  \varphi\left(  \mathbf{y}\right)
		\text{d}\mathbf{y\quad x}\in\Gamma_{j}, \quad j\in\left\{
		1,\cdots,q\right\};
	\end{align*}

	\item Double layer operators:%
	\begin{align*}
		K_{0}^{j,\ell}\psi\left(  \mathbf{x}\right)   &  :=\int_{\Gamma_{\ell}}%
		\frac{\partial G_{\kappa_{0}}\left(  \mathbf{x},\mathbf{y}\right)  }%
		{\partial\mathbf{n}_{0}\left(  \mathbf{y}\right)  }\varphi\left(
		\mathbf{y}\right)  \text{d}\mathbf{y\quad x}\in\Gamma_{j}\quad j,\ell \in\left\{
		0,\cdots,q\right\}  ,\\
		K_{j}\psi\left(  \mathbf{x}\right)   &  :=\int_{\Gamma_{j}}\frac{\partial
			G_{\kappa_{j}}\left(  \mathbf{x},\mathbf{y}\right)  }{\partial\mathbf{n}%
			_{j}\left(  \mathbf{y}\right)  }\varphi\left(  \mathbf{y}\right)
		\text{d}\mathbf{y\quad x}\in\Gamma_{j}, \quad j \in\left\{
		1,\cdots,q\right\};
	\end{align*}

	\item Adjoint double layer operators:%
	\begin{align*}
		\left(K^\prime\right)_{0}^{j,\ell}\varphi\left(  \mathbf{x}\right)   &  :=\frac{\partial
		}{\partial\mathbf{n}_{0}\left(  \mathbf{x}\right)  }\int_{\Gamma_{\ell}}%
		G_{\kappa_{0}}\left(  \mathbf{x},\mathbf{y}\right)  \varphi\left(
		\mathbf{y}\right)  \text{d}\mathbf{y\quad x}\in\Gamma_{j}\quad j,\ell\in\left\{
		0,\cdots,q\right\}  ,\\
		K_{j}^{\prime}\varphi\left(  \mathbf{x}\right)   &  :=\frac{\partial}%
		{\partial\mathbf{n}_{j}\left(  \mathbf{x}\right)  }\int_{\Gamma_{j}}%
		G_{\kappa_{j}}\left(  \mathbf{x},\mathbf{y}\right)  \varphi\left(
		\mathbf{y}\right)  \text{d}\mathbf{y\quad x}\in\Gamma_{j},\quad j \in\left\{
		1,\cdots,q\right\};
	\end{align*}

	\item Hypersingular operators:%
	\begin{align*}
		W_{0}^{j,\ell}\psi\left(  \mathbf{x}\right)   &  :=\frac{\partial}%
		{\partial\mathbf{n}_{0}\left(  \mathbf{x}\right)  }\int_{\Gamma_{\ell}}%
		\frac{\partial G_{\kappa_{0}}\left(  \mathbf{x},\mathbf{y}\right)  }%
		{\partial\mathbf{n}_{0}\left(  \mathbf{y}\right)  }\varphi\left(
		\mathbf{y}\right)  \text{d}\mathbf{y\quad x}\in\Gamma_{j}\quad j,\ell\in\left\{
		0,\cdots,q\right\}  ,\\
		W_{j}\psi\left(  \mathbf{x}\right)   &  :=\frac{\partial}{\partial
			\mathbf{n}_{j}\left(  \mathbf{x}\right)  }\int_{\Gamma_{j}}\frac{\partial
			G_{\kappa_{j}}\left(  \mathbf{x},\mathbf{y}\right)  }{\partial\mathbf{n}%
			_{j}\left(  \mathbf{y}\right)  }\varphi\left(  \mathbf{y}\right)
		\text{d}\mathbf{y\quad x}\in\Gamma_{j},\quad j \in\left\{
		1,\cdots,q\right\}.
	\end{align*}
\end{itemize}

{\color{black} To reformulate \eqref{transmission_problem_1}--\eqref{transmission_problem_4} in terms of a boundary integral equation, we consider the following Green's representations on the boundaries for the interior problems}: for $j = 1,\cdots,q$
\begin{subnumcases}{}
-\Vj \uNi(\x) - \frac{1}{2}\uDj(\x) + \Kj \uDj(\x) = & $0, \qquad \x\in\Gamma_j$ \label{BIE_1}
	\\
\frac{1}{2}\uNj(\x) + \Kpj\uNj(\x) - \Wj \uDj(\x) = & $0, \qquad \x\in\Gamma_j$, \label{BIE_2}
\end{subnumcases}
and for $j = 0,\cdots,q$

\begin{subnumcases}{}
	-\sum_{\ell = 0}^q\Vojl \uNol(\x) - \frac{1}{2}\uDoj(\x) + \sum_{\ell = 0}^q\Kojl \uDol(\x) = & $0, \qquad \x\in\Gamma_j$ \label{BIE_3}
	\\
	\frac{1}{2}\uNoj(\x) + \sum_{\ell = 0}^q\Kpojl\uNol(\x) - \sum_{\ell = 0}^q\Wojl \uDol(\x) = & $0, \qquad \x\in\Gamma_j$, \label{BIE_4}
\end{subnumcases}	
{\color{black} together with the transmission conditions \eqref{transmission_problem_2}--\eqref{transmission_problem_3}}.
Moreover, the boundary condition \eqref{transmission_problem_4} on $\Gammao$ reads {\color{black}(see, e.g. \cite{SauterSchwab2010}, Sec. 3.7)}
\begin{equation}\label{eq:DtN_map}
	-\Vooo \uNoo(\x) + \frac{1}{2}\uDoo(\x) + \Kooo \uDoo(\x) = \Vooo g(\x), \quad \x\in\Gammao.
\end{equation}
To write the variational formulation of the problem, reformulated in terms of \eqref{BIE_1}-\eqref{BIE_2}-\eqref{BIE_3}-\eqref{BIE_4}-\eqref{eq:DtN_map}, we introduce for a closed (not necessarily
connected) curve $\Gamma$ the space $X\left(  \Gamma\right)  =H^{1/2}\left(
\Gamma\right)  \times H^{-1/2}\left(  \Gamma\right)  $ and the scalar product
$\left\langle \cdot,\cdot\right\rangle _{X\left(  \Gamma\right)  }$ defined
as: for $\underline{w}=\left[  w_{D},w_{N}\right]  ^{T}$ and $\underline{v}%
=\left[  v_{D},v_{N}\right]  ^{T}\in X\left(  \Gamma\right)  $%
\[
\left\langle \underline{w},\underline{v}\right\rangle _{X\left(
	\Gamma\right)  }=\left\langle w_{D},v_{N}\right\rangle _{\Gamma}+\left\langle
v_{D},w_{N}\right\rangle _{\Gamma},
\]
where $\langle  f , g\rangle_{{\Gamma}} := \displaystyle{\int_{\Gamma} f \overline{g}}$ denotes
the $L^{2}$ inner product and its extension by density to the $H^{1/2}\left(
\Gamma\right)  \times H^{-1/2}\left(  \Gamma\right)  $ anti-duality pairing.

We first formulate the skeleton integral equations in the multi-trace space
${X^{\operatorname*{mult}}:=\prod_{j=0}^{q}X\left(  \partial\Omega_{j}\right)}$ and identify $X\left(  \partial\Omega_{0}\right)$ with the product space
${\prod_{j=0}^{q}X\left(  \Gamma_{j}\right)  }$ so that ${X^{\operatorname*{mult}%
}\cong\prod_{j=0}^{q}X\left(  \Gamma_{j}\right)  \times\prod_{j=1}^{q}X\left(
\Gamma_{j}\right)}$. We use the shorthands $X_{j}:=X\left(  \Gamma
_{j}\right)  $ and $\left\langle \underline{w},\underline{v}%
\right\rangle _{X_{j}} := \left\langle \underline{w},\underline{v}\right\rangle
_{X\left(  \Gamma_{j}\right)  }$.

Then, we introduce the operators
\begin{equation*}
\begin{aligned}
&C_j = \displaystyle{\begin{bmatrix}
	- \displaystyle{\frac{1}{2}{I} + \Kj} & \displaystyle{-\Vj}\\
	- \displaystyle{\Wj }&\displaystyle{\frac{1}{2}{I} + \Kpj}
	\end{bmatrix} }, & \qquad j = 1,\cdots,q, \\
&C_0^{j,\ell} = \displaystyle{\begin{bmatrix}
	\displaystyle{-\frac{1}{2}{\delta^{j,\ell}} + \Ko^{j,\ell}}  & \displaystyle{-\Vo^{j,\ell}}\\
	- \displaystyle{\Wo^{j,\ell}} & \displaystyle{\frac{1}{2}{\delta^{j,\ell}} +  \left(K^\prime\right)_0^{j,\ell}}
	\end{bmatrix} }, & \qquad j,\ell = 0,\cdots,q, \\
	&\widehat{C}_0^{0,0} = \displaystyle{\begin{bmatrix}
	\displaystyle{\frac{1}{2}{I} + \Ko^{0,0}}  & \displaystyle{-\Vo^{0,0}}\\
	\mathbf{0} & \mathbf{0}
	\end{bmatrix}} &
\end{aligned}
\end{equation*}
which allow {\color{black} us} to rewrite the variational forms of \eqref{BIE_1}-\eqref{BIE_2}, \eqref{BIE_3}-\eqref{BIE_4} and \eqref{eq:DtN_map}, respectively, in the compact form:
% NEW
the function $\left(  \underline{u}_{0}%
^{0},\underline{u}_{0}^{1},\ldots,\underline{u}_{0}^{q},\underline{u}%
_{1},\underline{u}_{2},\ldots,,\underline{u}_{q}\right)  \in
X^{\operatorname*{mult}}$ satisfy for $j\in\left\{  1,2,\ldots,q\right\}  $,
$m\in\left\{  0,1,\ldots,q\right\}  $%
\begin{equation}\label{eq:variational_tot}
\left\langle C_{j}\underline{u}_{j},\underline{v}_{j}\right\rangle _{X_{j}%
}=0,\quad\sum_{\ell=0}^{q}\left\langle C_{0}^{m,\ell}\underline{u}_{0}^{\ell
},\underline{v}_{0}^{m}\right\rangle _{X_{j}}=0,\quad\left\langle
\widehat{C}_{0}^{0,0}\underline{u}_{0}^{\ell},\underline{v}_{0}^{\ell
}\right\rangle _{X_{0}}=\left\langle \underline{g}_{0},\underline{v}_{0}%
^{0}\right\rangle _{X_{0}},
\end{equation}
%%%%%
where we have set $\ugo = [\Vooo g,0]^T$. {\color{black} We recall that in \eqref{eq:variational_tot} the transmission conditions \eqref{transmission_problem_2}--\eqref{transmission_problem_3} have also to be imposed}.

\noindent
By summing up the terms in \eqref{eq:variational_tot}, we aim at determining $\left(  \underline{u}_{0}%
^{0},\underline{u}_{0}^{1},\ldots,\underline{u}_{0}^{q},\underline{u}%
_{1},\underline{u}_{2},\ldots,,\underline{u}_{q}\right)  \in
X^{\operatorname*{mult}}$ such that
$$\underbrace{\sum_{j=0}^q \sum_{\ell = 0}^q\langle C_0^{j,\ell}\underline{u}_0^\ell,\uvo^j\rangle_{X_j} + \langle \widehat{C}_0^{0,0} \uuo^0, \uvo^0 \rangle_{X_0}}_{S_0} + \underbrace{\sum_{j=1}^q \langle C^j\underline{u}_j,\uvj\rangle_{X_j}}_{S_1}  = 
%S_0(\uuo,\uvo) + \sum_{i=1}^q S_i(\uui,\uvi) =  
\langle \underline{g}_0, \uvo\rangle_{X_0},$$
for all  $\left(\uvo^0,\uvo^1,\cdots,\uvo^q,\uvuno,\cdots,\uvq\right)\in
X^{\operatorname*{mult}}$, together with the transmission conditions \eqref{transmission_problem_2}--\eqref{transmission_problem_3}.
By making explicit the terms in $S_0$, we write

\begin{eqnarray}\label{eq:So_split}
S_0 &=& \sum_{j=0}^q \sum_{\ell = 0}^q\langle C_0^{j,\ell}\underline{u}_0^\ell,\uvo^j\rangle_{X_j} + \langle D_0^{0,0} \uuo^0, \uvo^0 \rangle_{X_0} \nonumber\\
&=& \sum_{j=0}^q\left[\langle - \frac{1}{2}\uDoj  + \sum_{\ell=0}^q\left(-\Vojl \uNol + \Kojl \uDol\right),\vNoj\rangle _{\Gamma_j} + \langle \vDoj,\frac{1}{2}\uNoj + \sum_{\ell=0}^q\left(\Kpojl\uNol - \Wojl \uDol\right)\rangle _{\Gamma_j}\right. \nonumber\\
  &+& \left.\langle \frac{1}{2}\uDoo + \Kooo\uDoo - \Vooo\uNoo, \vNoo\rangle _{\Gamma_0}  \right ] \nonumber\\
  &=& \sum_{j=0}^q\left[\langle - \frac{1}{2}\uDoj  - \Vojo \uNoo + \Kojo \uDoo,\vNoj\rangle _{\Gamma_j} +  \langle \sum_{\ell=1}^q\left(-\Vojl \uNol + \Kojl \uDol\right),\vNoj\rangle _{\Gamma_j}\right.\nonumber\\
 &+&  \left.\langle \vDoj,\frac{1}{2}\uNoj + \Kpojo\uNoo - \Wojo \uDoo\rangle _{\Gamma_j} + \langle \vDoj,\sum_{\ell=1}^q\left(\Kpojl\uNol - \Wojl \uDol\right)\rangle _{\Gamma_j}\right.\nonumber\\
  &+& \left.\langle \frac{1}{2}\uDoo + \Kooo\uDoo - \Vooo\uNoo, \vNoo\rangle _{\Gamma_0} \right ] \nonumber\\ 
 &=& \langle \cancel{- \frac{1}{2}\uDoo}  - \Vooo \uNoo + \Kooo \uDoo,\vNoo\rangle _{\Gamma_0}+  \langle \sum_{\ell=1}^q\left(-\Vool \uNol + \Kool \uDol\right),\vNoo\rangle _{\Gamma_0}\nonumber\\
 &+&  \langle \vDoo,\frac{1}{2}\uNoo + \Kpooo\uNoo - \Wooo \uDoo\rangle _{\Gamma_0}
 + \langle \vDoo,\sum_{\ell=1}^q\left(\Kpool\uNol - \Wool \uDol\right)\rangle _{\Gamma_0}\nonumber\\
 &+& \sum_{j=1}^q\left[\langle - \frac{1}{2}\uDoj  - \Vojo \uNoo + \Kojo \uDoo,\vNoj\rangle _{\Gamma_j}+  \langle \sum_{\ell=1}^q\left(-\Vojl \uNol + \Kojl \uDol\right),\vNoj\rangle _{\Gamma_j}\right.\nonumber\\
 &+& \left.\langle \vDoj,\frac{1}{2}\uNoj + \Kpojo\uNoo - \Wojo \uDoo\rangle _{\Gamma_j}+ \langle \vDoj,\sum_{\ell=1}^q\left(\Kpojl\uNol - \Wojl \uDol\right)\rangle _{\Gamma_j}\right]\nonumber\\
  &+& \left.\langle \cancel{\frac{1}{2}\uDoo} + \Kooo\uDoo - \Vooo\uNoo, \vNoo\rangle _{\Gamma_0} \right ].
\end{eqnarray}
We use the transmission relations \eqref{transmission_problem_2}-\eqref{transmission_problem_3} in 
\eqref{eq:So_split}, for both the test and trial functions, namely $\uDoj = \uDj$, $\vDoj = \vDj$, $\uNoj = -\uNj$ and $\vNoj = -\vNj$ for $j = 1,\cdots,q$. 
These relations allow us to formulate the skeleton integral equations in the
single-trace space $X^{\operatorname{single}}:=\prod_{j=0}^{q}X_{j}$ by
eliminating the one-sided traces on $\Gamma_{j}$ (from the side $\Omega_{0}$):

\begin{eqnarray}\label{eq:So_split2}
S_0
&=& \langle  - 2\Vooo \uNoo + 2\Kooo \uDoo,\vNoo\rangle _{\Gamma_0}+  \langle \sum_{\ell=1}^q\left(\Vool \uNl + \Kool \uDl\right),\vNoo\rangle _{\Gamma_0}\qquad\qquad\qquad\qquad\qquad\nonumber\\
&+&  \langle \vDoo,\frac{1}{2}\uNoo + \Kpooo\uNoo - \Wooo \uDoo\rangle _{\Gamma_0}
+ \langle \vDoo,\sum_{\ell=1}^q\left(-\Kpool\uNl - \Wool \uDl\right)\rangle _{\Gamma_0}\nonumber\\
&+& \sum_{j=1}^q\left[\langle - \frac{1}{2}\uDj  - \Vojo \uNoo + \Kojo \uDoo,-\vNj\rangle _{\Gamma_j}+  \langle \sum_{\ell=1}^q\left(\Vojl \uNl + \Kojl \uDl\right),-\vNj\rangle _{\Gamma_j}\right.\nonumber\\
&+& \left.\langle \vDj,-\frac{1}{2}\uNj + \Kpojo\uNoo - \Wojo \uDoo\rangle _{\Gamma_j}+ \langle \vDj,\sum_{\ell=1}^q\left(-\Kpojl\uNl - \Wojl \uDl\right)\rangle _{\Gamma_j}\right].
\end{eqnarray}
Then, summing up $S_0$ and $S_1$, we get
\begin{eqnarray}\label{eq:sum_S0_Sj}
S_0+ S_1 &=& \langle  - 2\Vooo \uNoo + 2\Kooo \uDoo,\vNoo\rangle _{\Gamma_0}+  \langle \sum_{\ell=1}^q\left(\Vool \uNl + \Kool \uDl\right),\vNoo\rangle _{\Gamma_0}\nonumber\\
&+&  \langle \vDoo,\frac{1}{2}\uNoo + \Kpooo\uNoo - \Wooo \uDoo\rangle _{\Gamma_0}
+ \langle \vDoo,\sum_{\ell=1}^q\left(-\Kpool\uNl - \Wool \uDl\right)\rangle _{\Gamma_0}\nonumber\\
&+& \sum_{j=1}^q\left[\langle \cancel{- \frac{1}{2}\uDj}  - \Vojo \uNoo + \Kojo \uDoo,-\vNj\rangle _{\Gamma_j}+  \langle \sum_{\ell=1}^q\left(\Vojl \uNl + \Kojl \uDl\right),-\vNj\rangle _{\Gamma_j}\right.\nonumber\\
&+& \left.\langle \vDj,\cancel{-\frac{1}{2}\uNj} + \Kpojo\uNoo - \Wojo \uDoo\rangle _{\Gamma_j}+ \langle \vDj,\sum_{\ell=1}^q\left(-\Kpojl\uNl - \Wojl \uDl\right)\rangle _{\Gamma_j}\right]\nonumber\\
&+&\sum_{j=1}^q \left[\langle -\Vj \uNj \cancel{- \frac{1}{2}\uDj} + \Kj \uDj,\vNj\rangle _{\Gamma_j} + \langle \vDj,\cancel{\frac{1}{2}\uNj} + \Kpj\uNj - \Wj \uDj\rangle _{\Gamma_j} \right]\nonumber\\
 & = &\langle \Vooo g,\vNoo\rangle _{\Gamma_0}.
\end{eqnarray}
Finally, the skeleton integral equation in the single-trace variational formulation reads: 

find $\uuu= [\uuo^0,\uuuno,\cdots,\uuq]\in X^{\operatorname{single}}$ such that 

\begin{eqnarray}\label{eq:sum_S0_Sj_2}
\uuC(\uuu,\uuv) &:=&  \sum_{j=1}^q\left[\langle -\Vojo \uNoo + \Kojo \uDoo,-\vNj\rangle _{\Gamma_j}+  \langle \sum_{\ell=1}^q\left(\Vojl \uNl + \Kojl \uDl\right),-\vNj\rangle _{\Gamma_j}\right.\nonumber\\
&+& \left.\langle \vDj, \Kpojo\uNoo - \Wojo \uDoo\rangle_{\Gamma_j}+ \langle \vDj,\sum_{\ell=1}^q\left(-\Kpojl\uNl - \Wojl \uDl\right)\rangle _{\Gamma_j}\right]\nonumber\\
&+&\sum_{j=1}^q \left[\langle -\Vj \uNj + \Kj \uDj,\vNj\rangle_{\Gamma_j} + \langle \vDj,\Kpj\uNj - \Wj \uDj\rangle _{\Gamma_j} \right]\nonumber\\
&+&\langle  - 2\Vooo \uNoo + 2\Kooo \uDoo,\vNoo\rangle _{\Gamma_0}+  \langle \sum_{\ell=1}^q\left(\Vool \uNl + \Kool \uDl\right),\vNoo\rangle _{\Gamma_0}\nonumber\\
&+&  \langle \vDoo,\frac{1}{2}\uNoo + \Kpooo\uNoo - \Wooo \uDoo\rangle _{\Gamma_0}
+ \langle \vDoo,\sum_{\ell=1}^q\left(-\Kpool\uNl - \Wool \uDl\right)\rangle _{\Gamma_0}\nonumber\\
 & = &\langle \Vooo g,\vNoo\rangle _{\Gamma_0}
\end{eqnarray}
for all $\uuv = [\uvo^0,\uvuno,\cdots,\uvq] \in X^{\operatorname{single}}$.
{\color{black}The well-posedness of problem \eqref{eq:sum_S0_Sj_2} follows, e.g., from \cite[Sec.
4]{Hiptmair_multiple_trace}, \cite[Sec. 4]{GrHiSa_pw_Lip}}.

% per operatori
\newcommand{\Vounoo}{V_{0}^{1,0}}
\newcommand{\Voouno}{V_{0}^{0,1}}
\newcommand{\Kounoo}{K_{0}^{1,0}}
\newcommand{\Koouno}{K_{0}^{0,1}}
\newcommand{\vNuno}{v_{N;1}}
\newcommand{\vDuno}{v_{D;1}}
\newcommand{\uNuno}{u_{N;1}}
\newcommand{\uDuno}{u_{D;1}}
\newcommand{\Vounouno}{V_{0}^{1,1}}
\newcommand{\Kounouno}{K_{0}^{1,1}}
\newcommand{\Kpounoo}{(K^\prime)_{0}^{1,0}}
\newcommand{\Kpoouno}{(K^\prime)_{0}^{0,1}}
\newcommand{\Kpounouno}{(K^\prime)_{0}^{1,1}}
\newcommand{\Wounoo}{W_{0}^{1,0}}
\newcommand{\Woouno}{W_{0}^{0,1}}
\newcommand{\Wounouno}{W_{0}^{1,1}}
\newcommand{\Vuno}{V_1}
\newcommand{\Kuno}{K_1}
\newcommand{\Kpuno}{K^\prime_1}
\newcommand{\Wuno}{W_1}
% per matrici e vettori
\newcommand{\vvNuno}{{\bf v}_{N;1}}
\newcommand{\vvDuno}{{\bf v}_{D;1}}
\newcommand{\uuNuno}{{\bf u}_{N;1}^M}
\newcommand{\uuDuno}{{\bf u}_{D;1}^M}
\newcommand{\uuDoo}{{\bf u}^{0,M}_{D;0}}
\newcommand{\uuNoo}{{\bf u}^{0,M}_{N;0}}
\newcommand{\VVounoo}{{\mathbb{V}_{0}^{1,0}}}
\newcommand{\VVoouno}{\mathbb{V}_{0}^{0,1}}
\newcommand{\KKounoo}{\mathbb{K}_{0}^{1,0}}
\newcommand{\KKoouno}{\mathbb{K}_{0}^{0,1}}
\newcommand{\VVounouno}{\mathbb{V}_{0}^{1,1}}
\newcommand{\KKounouno}{\mathbb{K}_{0}^{1,1}}
\newcommand{\KKpounoo}{(\mathbb{K}^\prime)_{0}^{1,0}}
\newcommand{\KKpoouno}{(\mathbb{K}^\prime)_{0}^{0,1}}
\newcommand{\KKpounouno}{(\mathbb{K}^\prime)_{0}^{1,1}}
\newcommand{\KKooo}{\mathbb{K}_{0}^{0,0}}
\newcommand{\KKpooo}{(\mathbb{K}^\prime)_{0}^{0,0}}
\newcommand{\VVooo}{\mathbb{V}_{0}^{0,0}}
\newcommand{\WWounoo}{\mathbb{W}_{0}^{1,0}}
\newcommand{\WWoouno}{\mathbb{W}_{0}^{0,1}}
\newcommand{\WWooo}{\mathbb{W}_{0}^{0,0}}
\newcommand{\WWounouno}{\mathbb{W}_{0}^{1,1}}
\newcommand{\VVuno}{\mathbb{V}_1}
\newcommand{\KKuno}{\mathbb{K}_1}
\newcommand{\KKpuno}{\mathbb{K}^\prime_1}
\newcommand{\WWuno}{\mathbb{W}_1}

\newcommand{\VVi}{\mathbb{V}_i}
\newcommand{\KKi}{\mathbb{K}_i}
\newcommand{\KKpi}{\mathbb{K}^\prime_i}
\newcommand{\WWi}{\mathbb{W}_i}
\newcommand{\Mi}{\mathcal{M}_i}
\newcommand{\Mo}{\mathcal{M}_0}
\newcommand{\MMi}{\mathbb{M}_i}
\newcommand{\MMo}{\mathbb{M}_0}

\newcommand{\bo}{\mathbf{0}}

\section{Concentric circles model problem}\label{sec:concentric_cirles}
In this section, we present the operator form of the skeleton integral equation for the specific case of a single transmission interface (which corresponds to the case $q=1$). Focusing on the model problem where $\Omega_1$ and $\Omega_0$ are concentric circles, we analyse the associated integral operators, derive the analytic solution representation, and perform an error analysis using a spectral Galerkin approximation method.

Let $0<r_{1}<r_{0}=1$. Let $\Omega_{1}$ be the open disc around the origin
with radius $r_{1}$ and boundary $\Gamma_{1}$ and let $\Omega$ be the open
disc around the origin with radius $r_{0}$ and boundary $\Gamma_{0}$. Then we set
$\Omega_{0}:=\Omega\backslash\overline{\Omega_{1}}$. 
In this case, the unknown $\uuu \in X^{\operatorname{single}}$ of \eqref{eq:sum_S0_Sj_2} consists of four terms $\uuu = [\uDoo,\uNoo,\uDuno,\uNuno]$ which are solution of the skeleton equations
\begin{eqnarray}\label{eq:skeleton_eq}
	&	\hskip-2cm\left[
	\begin{array}{cccc}
	-\Wounouno - \Wuno & \Kpuno -\Kpounouno &	-\Wounoo & \Kpounoo \\
		\\
		-\Kounouno+\Kuno & -\Vounouno-\Vuno &-\Kounoo & \Vounoo  \\
		\\
		-\Woouno & -\Kpoouno & -\Wooo & \frac{1}{2}{\mathbb D}+\Kpooo \\
		\\
		 \Koouno & \Voouno & 2\Kooo & -2\Vooo\\
	\end{array}
	\right]
	\, \left[
	\begin{array}{c}
		\uDuno\\
		\\
		\uNuno\\
		\\
		\uDoo\\
		\\
		\uNoo\\	
	\end{array}
	\right] 
	&=
	\left[
	\begin{array}{c}
		0 \\
		\\
		0\\
		\\
		0\\
		\\
		\Vooo g\\
	\end{array}
	\right].
\end{eqnarray}

\newcommand{\hum}{\hat{u}_{m}}
\newcommand{\huo}{\hat{u}_0}
\newcommand{\huuno}{\hat{u}_1}
\newcommand{\humi}{\hat{u}_{m,i}}
\newcommand{\humo}{\hat{u}_{m,0}}
\newcommand{\humuno}{\hat{u}_{m,1}}
\newcommand{\hui}{\hat{u}_{i}}
\newcommand{\Aiuno}{A_{i,1}}
\newcommand{\Aidue}{A_{i,2}}
\newcommand{\fmuno}{H^{(1)}_{m}}
\newcommand{\fmdue}{J_{m}}

In the following section we will study the properties of the integral operators involved in \eqref{eq:skeleton_eq}, using spectral approximating functions.
\subsection{Diagonal forms of skeleton integral operators}

\begin{lemma}
	For $m\in\mathbb{Z}$, let%
	\begin{equation}
		g_{m}\left(\x\right)  =\left(  \frac{x_{1}+\operatorname*{i}x_{2}}{\left\Vert
			\mbox{\boldmath$\x$}%
			\right\Vert }\right)  ^{m}\quad\text{for }%
		\mbox{\boldmath$\x$ }= (x_1,x_2)%
		%EndExpansion
		\in\mathbb{R}^{2}\backslash\left\{  \mathbf{0}\right\}  .\label{Defgm}%
	\end{equation}
	Let $\Gamma_{R}$ denote the circle with radius $R>0$ about the origin. Then%
	\begin{equation}\label{eq:(gm_gn)}
%	\left(  g_{m},g_{n}\right)  _{L^{2}\left(  \Gamma_{R}\right)  }=\delta
%	_{m,n}2\pi R.
{\color{black}\langle g_{m},g_{n}\rangle_{\Gamma_{R}}}=\delta
_{m,n}2\pi R.
	\end{equation}
	
\end{lemma}%

\proof

We use the parametrization: $\chi:\left[  -\pi,\pi\right[  \rightarrow
\Gamma_{R}$ given by $\chi\left(  t\right)  :=R\operatorname*{e}%
^{\operatorname*{i}t}$. Then%
\[
%\left(  g_{m},g_{n}\right)  _{L^{2}\left(  \Gamma_{R}\right)  }=R\int_{-\pi
%}^{\pi}\operatorname*{e}\nolimits^{\operatorname*{i}\left(  m-n\right)
%	t}dt=\delta_{m,n}2\pi R.
{\color{black}\langle g_{m},g_{n}\rangle_{\Gamma_{R}}}=R\int_{-\pi
}^{\pi}\operatorname*{e}\nolimits^{\operatorname*{i}\left(  m-n\right)
	t}dt=\delta_{m,n}2\pi R.
\]
\endproof

Let $J_m$ denote the Bessel function of order $m$ and $H^{(1)}_m$ denote the Hankel function of the first kind and order $m$ (see \cite{NIST:DLMF}). The following proposition is taken from \cite[Thm. 8]{Amini}, where the eigensystem of the integral operators is provided for the unit {\color{black}circle}. 

{\color{black}
\begin{proposition}
	\label{Propll}
	Let us define the
	functions%
	\[%
	\begin{array}
	[c]{ll}%
	\lambda_m^{V}\left(  x\right)  :=\displaystyle{\frac{\operatorname*{i}\pi}{2}J_{m}\left(
		x\right)  H_{m}^{\left(  1\right)  }\left(  x\right)}  , & \lambda_m^{K}\left(
	x\right)  :=-\displaystyle{\frac{\operatorname*{i}\pi}{2}xJ_{m}^{\prime}\left(  x\right)
		H_{m}^{\left(  1\right)  }\left(  x\right)} ,\\
	\lambda_m^{K^{\prime}}\left(  x\right)  :=\lambda_m^{K}\left(  x\right)  , &
	\lambda_m^{W}\left(  x\right)  :=\displaystyle{\frac{\operatorname*{i}\pi}{2}x^{2}%
		J_{m}^{\prime}\left(  x\right)  \left(  H_{m}^{\left(  1\right)  }\right)
		^{\prime}\left(  x\right) } .
	\end{array}
	\]
	
	Then, it holds%
	\[%
	\begin{array}
		[c]{ll}%
		\displaystyle{V_{0}^{\ell,\ell}g_{m}=r_{\ell}\lambda_{m}^{V}\left(  \kappa_{0}r_{\ell
		}\right)  g_{m}, \qquad \ell \in\{0,1\}}& \qquad \displaystyle{V_{1}g_{m}=r_{1}\lambda_{m}^{V}\left(  \kappa_{1}%
		r_{1}\right)  g_{m},}\\
		\\
		K_{0}^{1,1}g_{m}=\left(K^\prime\right)_{0}^{1,1}g_{m}=-\left(  \frac{1}{2}+\lambda_{m}%
		^{K}\left(  \kappa_{0}r_{1}\right)  \right)  g_{m} & \qquad K_{1}g_{m}=K_{1}^{\prime
		}g_{m}=\left(  \frac{1}{2}+\lambda_{m}^{K}\left(  \kappa_{1}r_{1}\right)
		\right)  g_{m}\\
		\\
		K_{0}^{0,0}g_{m}=\left(K^\prime\right)_{0}^{0,0}g_{m}=\left(  \frac{1}{2}+\lambda_{m}%
		^{K}\left(  \kappa_{0}r_{0}\right)  \right)  g_{m} & \\
		\\
		W_{0}^{\ell,\ell}g_{m}=\frac{1}{r_{\ell}}\lambda_{m}^{W}\left(  \kappa
		_{0}r_{\ell}\right)  g_{m}, \qquad \ell \in\{0,1\} & \qquad W_{1}g_{m}=\frac{1}{r_{1}}\lambda_{m}^{W}\left(
		\kappa_{1}r_{1}\right)  g_{m}.
	\end{array}
	\]
	
\end{proposition}}%

\proof
Following \cite[Thm. 8]{Amini}, we obtain the formulae%
\[%
\begin{array}
	[c]{ll}%
	\displaystyle{V_{0}^{0,0}g_{m}=\lambda_{m}^{V}\left(  \kappa_{0}r_{0}\right)  g_{m}}, &
	\displaystyle{K_{0}^{0,0}g_{m}=\left(  \frac{1}{2}+\lambda_{m}^{K}\left(  \kappa_{0}%
	r_{0}\right)  \right)  g_{m}},\\
	\displaystyle{\left(K^\prime\right)_{0}^{0,0}g_{m}=\left(  \frac{1}{2}+\lambda_{m}^{K^{\prime}}\left(
	\kappa_{0}r_{0}\right)  \right)  g_{m}}, & \displaystyle{W_{0}^{0,0}g_{m}=\lambda_{m}%
	^{W}\left(  \kappa_{0}r_{0}\right)  g_{m}}.
\end{array}
\]
A simple variable transformation then yields%
\begin{align*}
	V_{0}^{1,1}g_{m}\left(  \mathbf{x}\right)   &  =\int_{\Gamma_{1}}G_{\kappa
		_{0}}\left(  \mathbf{x},\mathbf{y}\right)  g_{m}\left(  \mathbf{y}\right)
	d\mathbf{y}\overset{\mathbf{\tilde{y}}=\mathbf{y}/r_{1};\mathbf{\tilde{x}%
		}=\mathbf{x}/r_{1}}{=}r_{1}\int_{\Gamma_{0}}G_{\kappa_{0}r_{1}}\left(
	\mathbf{\tilde{x}},\mathbf{\tilde{y}}\right)  g_{m}\left(  \mathbf{\tilde{y}%
	}\right)  d\mathbf{\tilde{y}}\\
	&  =r_{1}\lambda_{m}^{V}\left(  \kappa_{0}r_{1}\right)  g_{m}\left(
	\mathbf{\tilde{x}}\right)  =r_{1}\lambda_{m}^{V}\left(  \kappa_{0}%
	r_{1}\right)  g_{m}\left(  \mathbf{x}\right)  .
\end{align*}
For the operator $K$ we obtain%
\begin{align*}
	K_{0}^{1,1}g_{m}\left(  \mathbf{x}\right)   &  =\int_{\Gamma_{1}}%
	\frac{\partial G_{\kappa_{0}}\left(  \mathbf{x},\mathbf{y}\right)  }%
	{\partial\mathbf{n}_{0}\left(  \mathbf{y}\right)  }g_{m}\left(  \mathbf{y}%
	\right)  d\mathbf{y}\overset{\mathbf{\tilde{y}}=\mathbf{y}/r_{1}%
		;\mathbf{\tilde{x}}=\mathbf{x}/r_{1}}{=}\int_{\Gamma_{0}}\frac{\partial
		G_{\kappa_{0}r_{1}}\left(  \mathbf{\tilde{x}},\mathbf{\tilde{y}}\right)
	}{\partial\mathbf{n}_{0}\left(  \mathbf{\tilde{y}}\right)  }g_{m}\left(
	\mathbf{\tilde{y}}\right)  d\mathbf{\tilde{y}}\\
	&  =-\left(  \frac{1}{2}+\lambda_{m}^{K}\left(  \kappa_{0}r_{1}\right)
	\right)  g_{m}\left(  \mathbf{\tilde{x}}\right)  =-\left(  \frac{1}{2}%
	+\lambda_{m}^{K}\left(  \kappa_{0}r_{1}\right)  \right)  g_{m}\left(
	\mathbf{x}\right)  .
\end{align*}
For the adjoint $K^{\prime}$ we get, in a similar fashion,%
\begin{align*}
	(K^\prime)_{0}^{1,1}g_{m}\left(  \mathbf{x}\right)   &  =\left.  \left(
	\frac{\partial}{\partial\mathbf{n}_{0}\left(  \mathbf{\tilde{x}}\right)  }%
	\int_{\Gamma_{0}}G_{\kappa_{0}r_{1}}\left(  \mathbf{\tilde{x}},\mathbf{\tilde
		{y}}\right)  g_{m}\left(  \mathbf{\tilde{y}}\right)  d\mathbf{\tilde{y}%
	}\right)  \right\vert _{\mathbf{\tilde{x}}=\mathbf{x}/r_{1}}\\
	&  =-\left(  \frac{1}{2}+\lambda_{m}^{K^{\prime}}\left(  \kappa_{0}%
	r_{1}\right)  \right)  g_{m}\left(  \mathbf{\tilde{x}}\right)  =-\left(
	\frac{1}{2}+\lambda_{m}^{K^{\prime}}\left(  \kappa_{0}r_{1}\right)  \right)
	g_{m}\left(  \mathbf{x}\right)  .
\end{align*}
For the hypersingular operator it holds%
\begin{align*}
	W_{0}^{1,1}g_{m}\left(  \mathbf{x}\right)   &  =\frac{\partial}{\partial
		\mathbf{n}_{1}\left(  \mathbf{x}\right)  }\int_{\Gamma_{1}}\left\langle 
	\mathbf{n}_{1}\left(  \mathbf{y}\right)  ,\nabla_{\mathbf{y}}G_{\kappa_{0}%
	}\left(  \mathbf{x},\mathbf{y}\right)  \right\rangle g_{m}\left(
	\mathbf{y}\right)  d\mathbf{y}\\
	&  =\frac{1}{r_{1}}\frac{\partial}{\partial\mathbf{n}_{1}\left(
		\mathbf{\tilde{x}}\right)  }\int_{\Gamma_{0}}\left\langle  \mathbf{n}%
	_{1}\left(  \mathbf{\tilde{y}}\right)  ,\nabla_{\mathbf{\tilde{y}}}%
	G_{\kappa_{0}r_{1}}\left(  \mathbf{\tilde{x}},\mathbf{\tilde{y}}\right)
	\right\rangle g_{m}\left(  \mathbf{\tilde{y}}\right)  d\mathbf{\tilde{y}}\\
	&  =\frac{1}{r_{1}}\lambda_{m}^{W}\left(  \kappa_{0}r_{1}\right)  g_{m}\left(
	\mathbf{x}\right)  .
\end{align*}
The relations for the operators $V_{1}$, $K_{1}$, $K_{1}^{\prime}$, $W_{1}$
follow by taking into account that on $\Gamma_{1}$ it holds $\mathbf{n}%
_{0}=-\mathbf{n}_{1}$.%
\endproof

For the evaluation of the potentials on different surfaces, the following Lemma holds.

\begin{lemma}\label{lm:potentials_diag}
	It holds%

	\begin{align}\label{eq:pot_diff_surf}
\begin{array}
[c]{l}%
{\color{black}\langle g_{n},V_{0}^{0,1}g_{m}\rangle_{\Gamma_{0}}
}=\displaystyle{\delta_{m,n}2\pi r_{0}r_{1}\lambda_{m}^{V}\left(  \kappa_{0}r_{1}\right)
	\frac{H_{\left\vert m\right\vert }^{\left(  1\right)  }\left(  \kappa_{0}%
		r_{0}\right)  }{H_{\left\vert m\right\vert }^{\left(  1\right)  }\left(
		\kappa_{0}r_{1}\right)  }},\\
{\color{black}\langle g_{n},V_{0}^{1,0}g_{m}\rangle_{\Gamma_{1}}
}=\displaystyle{\delta_{m,n}2\pi r_{0}r_{1}\lambda_{m}^{V}\left(  \kappa_{0}r_{0}\right)
	\frac{J_{\left\vert m\right\vert }\left(  \kappa_{0}r_{1}\right)
	}{J_{\left\vert m\right\vert }\left(  \kappa_{0}r_{0}\right)  }},\\
{\color{black}\langle  g_{n},K_{0}^{0,1}g_{m}\rangle_{ \Gamma_{0}}
}=-\displaystyle{\delta_{m,n}2\pi r_{0}\lambda_{m}^{K}\left(  \kappa_{0}r_{1}\right)
	\frac{H_{\left\vert m\right\vert }^{\left(  1\right)  }\left(  \kappa_{0}%
		r_{0}\right)  }{H_{\left\vert m\right\vert }^{\left(  1\right)  }\left(
		\kappa_{0}r_{1}\right)  }},\\
{\color{black}\langle g_{n},K_{0}^{1,0}D_{0}^{0}g_{m}\rangle_{\Gamma
	_{1}}  }=\displaystyle{\delta_{m,n}2\pi r_{1}\left(  1+\lambda_{m}^{K}\left(
	\kappa_{0}r_{0}\right)  \right)  \frac{J_{\left\vert m\right\vert }\left(
		\kappa_{0}r_{1}\right)  }{J_{\left\vert m\right\vert }\left(  \kappa_{0}%
		r_{0}\right)  }}%
\end{array}
\end{align}
	and%
	\begin{align}\label{eq:pot_diff_surf2}		
\begin{array}
[c]{l}%
{\color{black}\langle g_{n},(K^\prime)_{0}^{0,1}g_{m}\rangle_{\Gamma
	_{0}} }=-\displaystyle{\delta_{m,n}2\pi r_{0}r_{1}\kappa_{0}\lambda_{m}^{V}\left(
	\kappa_{0}r_{1}\right)  \frac{\left(  H_{\left\vert m\right\vert }^{\left(
			1\right)  }\right)  ^{\prime}\left(  \kappa_{0}r_{0}\right)  }{H_{\left\vert
			m\right\vert }^{\left(  1\right)  }\left(  \kappa_{0}r_{1}\right)  }},\\
{\color{black}\langle g_{n},(K^\prime)_{0}^{1,0}g_{m}\rangle_{\Gamma
	_{1}} }=\displaystyle{\delta_{m,n}2\pi r_{0}r_{1}\kappa_{0}\lambda_{m}^{V}\left(
	\kappa_{0}r_{0}\right)  \frac{J_{\left\vert m\right\vert }^{\prime}\left(
		\kappa_{0}r_{1}\right)  }{J_{\left\vert m\right\vert }^{\left(  1\right)
		}\left(  \kappa_{0}r_{0}\right) } },\\
{\color{black}\langle g_{n},W_{0}^{0,1}g_{m}\rangle_{\Gamma_{0}}
}=\displaystyle{\delta_{m,n}2\pi r_{0}\kappa_{0}\lambda_{m}^{K}\left(  \kappa_{0}%
	r_{1}\right)  \frac{\left(  H_{\left\vert m\right\vert }^{\left(  1\right)
		}\right)  ^{\prime}\left(  \kappa_{0}r_{0}\right)  }{H_{\left\vert
			m\right\vert }^{\left(  1\right)  }\left(  \kappa_{0}r_{1}\right)  }},\\
{\color{black}\langle g_{n},W_{0}^{1,0}g_{m}\rangle_{\Gamma_{1}}
}=\displaystyle{\delta_{m,n}2\pi r_{1}\kappa_{0}\left(  1+\lambda_{m}^{K}\left(  \kappa
	_{0}r_{0}\right)  \right)  \frac{J_{\left\vert m\right\vert }^{\prime}\left(
		\kappa_{0}r_{1}\right)  }{J_{\left\vert m\right\vert }\left(  \kappa_{0}%
		r_{0}\right) } }.
\end{array}
\end{align}				
			\end{lemma}%
			
			\proof
			{\color{black}Recalling the definition of the trace operators \eqref{eq:trace_D} and \eqref{eq:trace_N},  we make use of the well-known trace relations (see, e.g., \cite[Chap.
				3.3]{SauterSchwab2010}):%
\begin{equation*}
\begin{aligned}
&\gamma_{D;0}^0 S_{0}^{0}=V_{0}^{0,0}, && \qquad\qquad
\gamma_{D;0}^1 S_{0}^{1}=V_{0}^{1,1},\\
&\gamma_{D;1} S_{1}^{1}=V_{1}, & &\\
&\gamma_{D;0}^0 D_{0}^{0}=\frac{1}{2}%
I+K_{0}^{0,0}, & &\qquad \qquad\gamma_{D;0}^1 D_{0}%
^{1}=\frac{1}{2}I+K_{0}^{1,1},\\
&\gamma_{D;1} D_{1}^{1}=\frac{1}{2}%
I+K_{1}, && \\
&\gamma_{N;0}^0 S_{0}^{0}=-\left(  \frac
{1}{2}I+(K^\prime)_{0}^{0,0}\right),& & \qquad\qquad \gamma_{N;0}^1 S_{0}^{1}=-\left(  \frac{1}{2}I+(K^\prime)_{0}^{1,1
}\right), \\
&\gamma_{N;1} S_{1}^{1}=-\left(  \frac
{1}{2}I+K_{1}^{\prime}\right). & &
\end{aligned}
\end{equation*}
}
			Next, we consider the exterior Helmholtz Dirichlet problem%
			\begin{equation}%
				\begin{array}
					[c]{cl}%
					-\Delta u_{m}(\x)-\kappa_{0}^{2}u_{m}(\x) = 0 & \x\in\Omega^{e}:=\mathbb{R}%
					^{2}\backslash\overline{\Omega_{1}},\\
					%\Omega_{1}^{c}:=\mathbb{R}^{2}\backslash\overline{\Omega_{1}},\\
					u_{m}(\x) = g_{m}(\x) & \x\in\Gamma_{1},%\\
					%u_{m}\text{ satisfies SRC} & \text{at infinity.}%
				\end{array}
				\label{freHelm}%
			\end{equation}
			{\color{black}equipped with the \emph{Sommerfeld radiation condition}
			$$\lim\limits_{\|\x\|\rightarrow 0} \|\x\|^\frac{1}{2}\big(\langle \x / \|\x\|,\nabla u_m\rangle - \operatorname*{i}\kappa_0 u_m(\x)\big) = 0.$$}
			For the solution we employ a single layer ansatz $u_{m} = S_{0}^{1}\varphi_{m}$, 
			for a boundary density $\varphi_{m}$ which is determined by the boundary
			conditions, that is $V_{0}^{1,1}\varphi_{m} = g_{m}$.
			In this way, from Proposition \ref{Propll} we get%
			\[
			\varphi_{m}:=\left(  r_{1}\lambda_{m}^{V}\left(  \kappa_{0}r_{1}\right)
			\right)  ^{-1}g_{m}%
			\]
			and%
			\[
			u_{m}=S_{0}^{1}\varphi_{m}=\left(  r_{1}\lambda_{m}^{V}\left(  \kappa_{0}%
			r_{1}\right)  \right)  ^{-1}S_{0}^{1}g_{m}.
			\]
			On the other hand, it is known (see, e.g., \cite[(2.10)]{MonkChandlerWilde},
			\cite[Chap. 2.6]{Nedelec01}) that the solution of (\ref{freHelm}) is given by%
			\[
			u_{m}=\frac{H_{\left\vert m\right\vert }^{\left(  1\right)  }\left(
				\kappa_{0}r\right)  }{H_{\left\vert m\right\vert }^{\left(  1\right)  }\left(
				\kappa_{0}r_{1}\right)  }g_{m}%
			\]
			and a comparison leads to%
			\[
			S_{0}^{1}g_{m}=r_{1}\lambda_{m}^{V}\left(  \kappa_{0}r_{1}\right)
			\frac{H_{\left\vert m\right\vert }^{\left(  1\right)  }\left(  \kappa
				_{0}r\right)  }{H_{\left\vert m\right\vert }^{\left(  1\right)  }\left(
				\kappa_{0}r_{1}\right)  }g_{m}\quad\text{in }\Omega^e.
			\]			
			In a similar way, for the interior Helmholtz Dirichlet problems associated with $\Omega$ and $\Omega_1$, we compute, respectively%
			\begin{align*}
				S_{0}^{0}g_{m} &  =r_{0}\lambda_{m}^{V}\left(  \kappa_{0}r_{0}\right)
				\frac{J_{\left\vert m\right\vert }\left(  \kappa_{0}r\right)  }{J_{\left\vert
						m\right\vert }\left(  \kappa_{0}r_{0}\right)  }g_{m}\quad\text{in }\Omega,\\
				S_{1}^{1}g_{m} &  =r_{1}\lambda_{m}^{V}\left(  \kappa_{1}r_{1}\right)
				\frac{J_{\left\vert m\right\vert }\left(  \kappa_{1}r\right)  }{J_{\left\vert
						m\right\vert }\left(  \kappa_{1}r_{1}\right)  }g_{m}\quad\text{in }\Omega_{1}.
			\end{align*}
			Next, we compute $D_{0}^{\ell}g_{m}$, for $\ell = 0,1$ and start with an ansatz via the double
			layer potential:%
			\[
			u_{m}=D_{0}^{1}\psi_{m}.
			\]
			The boundary condition at $\Gamma_{1}$ leads to%
			\[
			\left(  \frac{1}{2}I+K_{0}^{1,1}\right)  \psi_{m}=g_{m}
			\]
		so that, from Proposition \ref{Propll}, we get $\psi_{m}=-\left(  \lambda_{m}^{K}\left(  \kappa_{0}r_{1}\right)  \right)
			^{-1}g_{m}$.
			In a similar fashion, as for the single layer potential, we obtain%
			\[
			D_{0}^{1}g_{m}=-\lambda_{m}^{K}\left(  \kappa_{0}r_{1}\right)  \frac
			{H_{\left\vert m\right\vert }^{\left(  1\right)  }\left(  \kappa_{0}r\right)
			}{H_{\left\vert m\right\vert }^{\left(  1\right)  }\left(  \kappa_{0}%
				r_{1}\right)  }g_{m}\quad\text{in }\Omega^e.
			\]
			The application of the same arguments to the remaining operators leads to%
			\begin{align*}
				D_{0}^{0}g_{m} &  =\left(  1+\lambda_{m}^{K}\left(  \kappa_{0}r_{0}\right)
				\right)  \frac{J_{\left\vert m\right\vert }\left(  \kappa_{0}r\right)
				}{J_{\left\vert m\right\vert }\left(  \kappa_{0}r_{0}\right)  }g_{m}%
				\quad\text{in }\Omega,\\
				D_{1}^{1}g_{m} &  =\left(  1+\lambda_{m}^{K}\left(  \kappa_{1}r_{1}\right)
				\right)  \frac{J_{\left\vert m\right\vert }\left(  \kappa_{1}r\right)
				}{J_{\left\vert m\right\vert }\left(  \kappa_{1}r_{1}\right)  }g_{m}%
				\quad\text{in }\Omega_{1}.
			\end{align*}
			From the above relations and recalling \eqref{eq:(gm_gn)}, relations \eqref{eq:pot_diff_surf} are obtained by applying the Dirichlet traces of the potentials:%
			
{\color{black}		
		\begin{align*}
	\langle  g_{n},\gamma_{D;0}^0 S_{0}^{1}g_{m}\rangle_{\Gamma_{0}  } &
	=\delta_{m,n}2\pi r_{0}r_{1}\lambda_{m}^{V}\left(  \kappa_{0}r_{1}\right)
	\frac{H_{\left\vert m\right\vert }^{\left(  1\right)  }\left(  \kappa_{0}%
		r_{0}\right)  }{H_{\left\vert m\right\vert }^{\left(  1\right)  }\left(
		\kappa_{0}r_{1}\right)  },\\
	\langle   g_{n},\gamma_{D;0}^1  S_{0}^{0}g_{m}\rangle_{\Gamma_{1} } &
	=\delta_{m,n}2\pi r_{0}r_{1}\lambda_{m}^{V}\left(  \kappa_{0}r_{0}\right)
	\frac{J_{\left\vert m\right\vert }\left(  \kappa_{0}r_{1}\right)
	}{J_{\left\vert m\right\vert }\left(  \kappa_{0}r_{0}\right)  },\\
	\langle   g_{n},\gamma_{D;0}^0 D_{0}^{1}g_{m}\rangle_{  \Gamma_{0} } &
	=-\delta_{m,n}2\pi r_{0}\lambda_{m}^{K}\left(  \kappa_{0}r_{1}\right)
	\frac{H_{\left\vert m\right\vert }^{\left(  1\right)  }\left(  \kappa_{0}%
		r_{0}\right)  }{H_{\left\vert m\right\vert }^{\left(  1\right)  }\left(
		\kappa_{0}r_{1}\right)  },\\
	\langle   g_{n},\gamma_{D;0}^1 D_{0}^{0}g_{m}\rangle_{ \Gamma_{1}  } &
	=\delta_{m,n}2\pi r_{1}\left(  1+\lambda_{m}^{K}\left(  \kappa_{0}%
	r_{0}\right)  \right)  \frac{J_{\left\vert m\right\vert }\left(  \kappa
		_{0}r_{1}\right)  }{J_{\left\vert m\right\vert }\left(  \kappa_{0}%
		r_{0}\right)  }.
	\end{align*}
}			
Finally, for the normal traces of the potentials we get%
			{\color{black}
					\begin{align*}
				\langle g_{n},\gamma_{N;0}^0 S_{0}^{1}%
				g_{m}\rangle  _{\Gamma_{0}  } &  =-\delta_{m,n}2\pi
				r_{0}r_{1}\kappa_{0}\lambda_{m}^{V}\left(  \kappa_{0}r_{1}\right)
				\frac{\left(  H_{\left\vert m\right\vert }^{\left(  1\right)  }\right)
					^{\prime}\left(  \kappa_{0}r_{0}\right)  }{H_{\left\vert m\right\vert
					}^{\left(  1\right)  }\left(  \kappa_{0}r_{1}\right)  },\\
				\langle g_{n},\gamma_{N;0}^1 S_{0}^{0}%
				g_{m}\rangle  _{\Gamma_{1}  } &  =\delta_{m,n}2\pi
				r_{0}r_{1}\kappa_{0}\lambda_{m}^{V}\left(  \kappa_{0}r_{0}\right)
				\frac{J_{\left\vert m\right\vert }^{\prime}\left(  \kappa_{0}r_{1}\right)
				}{J_{\left\vert m\right\vert }^{\left(  1\right)  }\left(  \kappa_{0}%
					r_{0}\right)  },\\
				\langle g_{n},\gamma_{N;0}^0 D_{0}^{1}%
				g_{m}\rangle  _{ \Gamma_{0}  } &  =\delta_{m,n}2\pi
				r_{0}\kappa_{0}\lambda_{m}^{K}\left(  \kappa_{0}r_{1}\right)  \frac{\left(
					H_{\left\vert m\right\vert }^{\left(  1\right)  }\right)  ^{\prime}\left(
					\kappa_{0}r_{0}\right)  }{H_{\left\vert m\right\vert }^{\left(  1\right)
					}\left(  \kappa_{0}r_{1}\right)  },\\
				\langle g_{n},\gamma_{N;0}^1 D_{0}^{0}%
				g_{m}\rangle  _{\Gamma_{1} } &  =\delta_{m,n}2\pi
				r_{1}\kappa_{0}\left(  1+\lambda_{m}^{K}\left(  \kappa_{0}r_{0}\right)
				\right)  \frac{J_{\left\vert m\right\vert }^{\prime}\left(  \kappa_{0}%
					r_{1}\right)  }{J_{\left\vert m\right\vert }\left(  \kappa_{0}r_{0}\right)  },
				\end{align*}%
			}
			from which \eqref{eq:pot_diff_surf2} follow. 
			\endproof

			\subsection{Regularity theory and spectral approximability}% for the circular transmission problem}
			
			\subsubsection{Representation of the traces of the analytic solution}
			
			In this section we analyze the regularity of the analytic solution for the
			transmission problem \eqref{transmission_problem_1}-\eqref{transmission_problem_4} in the case $\Omega_0$ and $\Omega_1$ are the domains introduced in Section \ref{sec:concentric_cirles}, and 
			for a given right-hand side $g\in H^{-1/2}\left(  \Gamma_{0}\right)  $. 
			In what follows, we use the convention that if $\mathbf{x}\in\mathbb{R}^{2}$ and $r\geq0$
			appear in the same context the relation $r=\left\Vert \mathbf{x}\right\Vert $ holds.
			We employ a Fourier analysis and first analyze the solution $u_{m}$ for the
			right-hand side $g_{m}$ as in (\ref{Defgm}). To obtain a representation
			of $u_{m}$, we introduce the coefficients%
			\[
			A_{m}^{1,1}:=0\quad\text{and\quad}A_{m}^{0,2}:=\frac{\pi\operatorname*{i}%
				r_{0}}{2}H_{m}^{\left(  1\right)  }\left(  \kappa_{0}r_{0}\right)  .
			\]
			The coefficients {\color{black}$A_{m}^{0,1}$ and $A_{m}^{1,2}$} are the solution of the linear system%
			\begin{equation}
				\left[
				\begin{array}
					[c]{cc}%
					H_{m}^{\left(  1\right)  }\left(  \kappa_{0}r_{1}\right)   & -J_{m}\left(
					\kappa_{1}r_{1}\right)  \\
					\kappa_{0}\left(  H_{m}^{\left(  1\right)  }\right)  ^{\prime}\left(
					\kappa_{0}r_{1}\right)   & -\kappa_{1}J_{m}^{\prime}\left(  \kappa_{1}%
					r_{1}\right)
				\end{array}
				\right]  \left(
				\begin{array}
					[c]{c}%
					A_{m}^{0,1}\\
					A_{m}^{1,2}%
				\end{array}
				\right)  =-A_{m}^{0,2}\left(
				\begin{array}
					[c]{l}%
					J_{m}\left(  \kappa_{0}r_{1}\right)  \\
					\kappa_{0}J_{m}^{\prime}\left(  \kappa_{0}r_{1}\right)
				\end{array}
				\right)  \label{SLP}%
			\end{equation}
			explicity given by%
			\begin{equation}
				\left(
				\begin{array}
					[c]{c}%
					A_{m}^{0,1}\\
					A_{m}^{1,2}%
				\end{array}
				\right)  =\frac{A_{m}^{0,2}}{\mathcal{W}\left(  H_{m}^{\left(  1\right)
					}\left(  \kappa_{0}\cdot\right),  J_{m}\left(  \kappa_{1}\cdot\right)  \right)
					\left(  r_{1}\right)  }\left(
				\begin{array}
					[c]{c}%
					\mathcal{W}\left(  J_{m}\left(  \kappa_{1}\cdot\right)  ,J_{m}\left(
					\kappa_{0}\cdot\right)  \right)  \left(  r_{1}\right)  \\
					-\frac{2\operatorname*{i}}{\pi r_{1}}%
				\end{array}
				\right)  ,\label{coeff0112}%
			\end{equation}
			where $\mathcal{W}\left(  f_{1},f_{2}\right)  \left(  z\right)  =f_{1}\left(
			z\right)  f_{2}^{\prime}\left(  z\right)  -f_{1}^{\prime}\left(  z\right)
			f_{2}\left(  z\right)  $ denotes the Wronskian of the two functions $f_{1}$,
			$f_{2}$ evaluated at $z$. 
			
			\begin{proposition}
				For $0<r_{1}<r_{0}=1$ and $\kappa_{0}$, $\kappa_{1}>0$ the linear system
				(\ref{SLP}) has a unique solution.
			\end{proposition}
			
			The proof can be found, e.g., in \cite[Lemma 1]{Poignard_vogelius}.
			
			\
			
			\noindent
			Note that (cf. (\cite[10.5.3]{NIST:DLMF})%
			\[
			\mathcal{W}\left(  H_{m}^{\left(  1\right)  }\left(  x\cdot\right),
			J_{m}\left(  x\cdot\right)  \right)  \left(  z\right)  =\frac{2}%
			{\pi\operatorname*{i}z}.
			\]
			This allows us to define the radial depending function $\hat{u}_{m}:\left[
			0,r_{0}\right]  \rightarrow\mathbb{C}$ by%
			\[
			\left.  \hat{u}_{m}\right\vert _{\tau_{j}}\left(  r\right)  :=\hat{u}%
			_{m,j}\left(  r\right)  :=A_{m}^{j,1}H_{m}^{\left(  1\right)  }\left(
			\kappa_{j}r\right)  +A_{m}^{j,2}J_{m}\left(  \kappa_{j}r\right)
			\]
			for $\tau_{1}:=\left[  0,r_{1}\right[  $, $\tau_{0}:=\left]  r_{1}%
			,r_{0}\right[  $ and finally the solution%
			\[
			u_{m}\left(  \mathbf{x}\right)  =\hat{u}_{m}\left(  r\right)  g_{m}(\x).
			\]

			\begin{remark}
				The coefficients $A_{m}^{0,1}$, $A_{m}^{1,2}$ can be expressed by the functions%
				\begin{align}
					R_{m}\left(  \kappa_{0}r_{1},\kappa_{1}/\kappa_{0}\right)   &  :=\frac
					{\mathcal{W}\left(  J_{m}\left(  \kappa_{1}\cdot\right)  ,J_{m}\left(
						\kappa_{0}\cdot\right)  \right)  \left(  r_{1}\right)  }{\mathcal{W}\left(
						H_{m}^{\left(  1\right)  }\left(  \kappa_{0}\cdot\right),  J_{m}\left(
						\kappa_{1}\cdot\right)  \right)  \left(  r_{1}\right)  },\nonumber\\
					T_{m}\left(  \kappa_{0}r_{1},\kappa_{1}/\kappa_{0}\right)   &  :=-\frac
					{2\operatorname*{i}\kappa_{0}}{\pi}\frac{1}{\mathcal{W}\left(  H_{m}^{\left(
							1\right)  }\left(  \kappa_{0}\cdot\right),  J_{m}\left(  \kappa_{1}%
						\cdot\right)  \right)  \left(  r_{1}\right)  },\nonumber\\
					S_{m}\left(  \kappa_{0}r_{1}\right)   &  :=-\frac{R_{m}\left(  \kappa_{0}%
						r_{1},\kappa_{1}/\kappa_{0}\right)  }{J_{m}\left(  \kappa_{0}r_{1}\right)
					}H_{m}^{\left(  1\right)  }\left(  \kappa_{0}r_{1}\right)  \label{Srel}%
				\end{align}
				as in \cite[(7.1)]{Capdebosq_ball} according to%
				\begin{align}
					\left(
					\begin{array}
						[c]{c}%
						A_{m}^{0,1}\\
						A_{m}^{1,2}%
					\end{array}
					\right)   &  =\frac{\pi\operatorname*{i}r_{0}}{2}\left(
					\begin{array}
						[c]{c}%
						R_{m}\left(  \kappa_{0}r_{1},\kappa_{1}/\kappa_{0}\right)  \\
						T_{m}\left(  \kappa_{0}r_{1},\kappa_{1}/\kappa_{0}\right)
					\end{array}
					\right)  H_{m}^{\left(  1\right)  }\left(  \kappa_{0}r_{0}\right)
					,\label{Am00}\\
					\left(
					\begin{array}
						[c]{c}%
						A_{m}^{0,1}H_{m}^{\left(  1\right)  }\left(  \kappa_{0}r_{1}\right)  \\
						A_{m}^{1,2}J_{m}\left(  \kappa_{1}r_{1}\right)  \\
						\kappa_{1}A_{m}^{1,2}J_{m}^{\prime}\left(  \kappa_{1}r_{1}\right)
					\end{array}
					\right)   &  =\frac{\pi\operatorname*{i}r_{0}}{2}\left(
					\begin{array}
						[c]{c}%
						-S_{m}\left(  \kappa_{0}r_{1}\right)  \\
						1-S_{m}\left(  \kappa_{0}r_{1}\right)  \\
						\kappa_{0}\left(  \frac{J_{m}^{\prime}\left(  \kappa_{0}r_{1}\right)  }%
						{J_{m}\left(  \kappa_{0}r_{1}\right)  }-S_{m}\left(  \kappa_{0}r_{1}\right)
						\frac{\left(  H_{m}^{\left(  1\right)  }\right)  ^{\prime}\left(  \kappa
							_{0}r_{1}\right)  }{H_{m}^{\left(  1\right)  }\left(  \kappa_{0}r_{1}\right)
						}\right)
					\end{array}
					\right)  J_{m}\left(  \kappa_{0}r_{1}\right)  H_{m}^{\left(  1\right)
					}\left(  \kappa_{0}r_{0}\right)  .\label{Am01}%
				\end{align}
				The functions $R_{m}$, $T_{m}$, $S_{m}$ are studied in detail in
				\cite{Capdebosq_ball}.
			\end{remark}
			For a given Dirichlet datum%
			\begin{equation}
				g(\x)=\sum_{m\in\mathbb{Z}}\alpha_{m}g_{m}(\x)\label{grep}%
			\end{equation}
			the corresponding solution is given by%
			\begin{equation}\label{u_ex}			
			u\left(  \mathbf{x}\right)  =\sum_{m\in\mathbb{Z}}\alpha_{m}\hat{u}_{m}\left(
			r\right)  g_{m}\left(  \mathbf{x}\right)  .
			\end{equation}
			Consequently the four traces which appear in our skeleton integral equation
			\eqref{eq:sum_S0_Sj_2} are given by%
						\begin{subequations}
							\label{userie}
						\end{subequations}%
						%EndExpansion%
						\begin{align}
							u_{\operatorname*{D};1} &  :=\sum_{m\in\mathbb{Z}}\alpha_{m}A_{m}^{1,2}%
							J_{m}\left(  \kappa_{1}r_{1}\right)  g_{m},\tag{%
								%TCIMACRO{\TeXButton{userie}{\ref{userie}}}%
								%BeginExpansion
								\ref{userie}%
								%EndExpansion
								a}\label{useriea}\\
							u_{\operatorname*{N};1} &  :=-\kappa_{1}\sum_{m\in\mathbb{Z}}\alpha_{m}%
							A_{m}^{1,2}J_{m}^{\prime}\left(  \kappa_{1}r_{1}\right)  g_{m},\tag{%
								%TCIMACRO{\TeXButton{userie}{\ref{userie}}}%
								%BeginExpansion
								\ref{userie}%
								%EndExpansion
								b}\label{userieb}\\
							u_{\operatorname*{D};0}^{0} &  :=\left.  u_{\operatorname*{D};0}\right\vert
							_{\Gamma_{0}}:=\sum_{m\in\mathbb{Z}}\alpha_{m}\left(  A_{m}^{0,1}%
							H_{m}^{\left(  1\right)  }\left(  \kappa_{0}r_{0}\right)  +A_{m}^{0,2}%
							J_{m}\left(  \kappa_{0}r_{0}\right)  \right)  g_{m},\tag{%
								%TCIMACRO{\TeXButton{userie}{\ref{userie}}}%
								%BeginExpansion
								\ref{userie}%
								%EndExpansion
								c}\label{useriec}\\
							u_{\operatorname*{N};0}^{0} &  :=\left.  u_{\operatorname*{N};0}\right\vert
							_{\Gamma_{0}}:=-\kappa_{0}\sum_{m\in\mathbb{Z}}\alpha_{m}\left(  A_{m}%
							^{0,1}\left(  H_{m}^{\left(  1\right)  }\right)  ^{\prime}\left(  \kappa
							_{0}r_{0}\right)  +A_{m}^{0,2}J_{m}^{\prime}\left(  \kappa_{0}r_{0}\right)
							\right)  g_{m}.\tag{%
								%TCIMACRO{\TeXButton{userie}{\ref{userie}}}%
								%BeginExpansion
								\ref{userie}%
								%EndExpansion
								d}\label{useried}%
						\end{align}						
{\color{black}Recall that $\kappa_{j}=\kappa\sqrt{n_{j}}$; if we consider the material
parameters $n_{j}$ and $r_{1}\in\left]  0,1\right[  $ as being fixed, the
\textit{resonance frequencies} $\kappa$ of the boundary value problem
\eqref{transmission_problem_1}--\eqref{transmission_problem_4} are the complex zeroes of
\begin{equation}\label{eq:wronskian}
\mathcal{W}\left(  H_{m}^{\left(  1\right)}\left(  \kappa_{0}\cdot\right),  J_{m}\left(  \kappa_{1}\cdot\right)  \right)\left(  r_{1}\right)  = -\kappa_1 \fmuno(\kappa_0 r_1)\fmdue^{\prime}(\kappa_1 r_1) + \kappa_0 \fmdue(\kappa_1 r_1)(\fmuno)^{\prime}(\kappa_0 r_1),
\end{equation}
as a function of the wavenumber $\kappa$. It is known from, e.g., \cite[Prop.
1]{Bouchra_Sa_Whis_I} (see also \cite{Poignard_vogelius}, \cite{Capdeboscq_3D}%
) that these zeroes satisfy $\operatorname{Im}\kappa<0$. Hence, for the
considered range of wave numbers $\kappa\in\mathbb{C}\backslash\left\{
0\right\}  $ with $\mathfrak{\operatorname{Im}}\left(  \kappa\right)  \geq0$,
the coefficients $A_{m}^{0,1}$, $A_{m}^{1,2}$ are well-defined by
\eqref{coeff0112}.}

\subsubsection{Estimates of the Fourier coefficients}
						
In this section we estimate the coefficients in the Fourier representations
(\ref{userie}) of the traces and normal derivatives of the solution. We start by observing that from
\cite[10.4.1, .2]{NIST:DLMF} it follows%
\begin{equation}%
\begin{array}
[c]{ll}%
J_{-m}\left(  z\right)  =\left(  -1\right)  ^{m}J_{m}\left(  z\right)  , &
H_{-m}^{\left(  1\right)  }\left(  z\right)  =\left(  -1\right)  ^{m}%
H_{m}^{\left(  1\right)  }\left(  z\right),
\end{array}
\label{reflBessel}%
\end{equation}
so that the Fourier coefficients $A_{m}^{i,2}J_{m}\left(  \kappa_{i}%
r_{i}\right)  $, $A_{m}^{i,2}J_{m}^{\prime}\left(  \kappa_{i}r_{i}\right)$, for $i=0,1$, $A_{m}^{0,1}H^{(1)}_{m}\left(  \kappa_{0}%
r_{0}\right)$ and $A_{m}^{0,1}\left(H_{m}^{(1)}\right)^{\prime}\left(  \kappa_{0}r_{0}\right)$
 in (\ref{userie}) only depend on $\left\vert m\right\vert $.

\begin{lemma}\label{lm:estimates_four_coeff}
\label{LemFC}For any $\eta\in\left]  0,1\right[  $, there exists a constant
$C_{\eta}$ such that for any $\mu>\operatorname*{e}/\left(  2\eta\right)  $ it holds:

\begin{enumerate}
\item[a.] for $\left\vert m\right\vert \geq\mu\kappa_{0}r_{0}>0$ it holds%
\[
\left\vert A_{m}^{0,2}J_{m}\left(  \kappa_{0}r_{0}\right)  \right\vert \leq
C_{\eta}\frac{r_{0}}{\left\vert m\right\vert };
\]

\item[b.] for $\left\vert m\right\vert \geq r_{1}\max\left\{  \mu\kappa
_{0},\kappa_{1}\right\}  $ it holds%
\[
\left\vert A_{m}^{0,1}H_{m}^{\left(  1\right)  }\left(  \kappa_{0}%
r_{0}\right)  \right\vert \leq C_{\eta}r_{0}r_{1}\frac{\left\vert \kappa
	_{0}-\kappa_{1}\right\vert }{\left\vert m\right\vert ^{4/3}}\operatorname*{e}%
\nolimits^{\kappa_{0}\left(  r_{0}-r_{1}\right)  }\left(  \frac{r_{1}}{r_{0}%
}\right)  ^{|m|};
\]

\item[c.] for $\left\vert m\right\vert \geq2+r_{1}\max\left\{  \mu\kappa
_{0},\kappa_{1}\right\}  $ it holds%
\[
\left\vert \kappa_{0}A_{m}^{0,1}\left(  H_{m}^{\left(  1\right)  }\right)
^{\prime}\left(  \kappa_{0}r_{0}\right)  \right\vert \leq C_{\eta}r_{1}%
\frac{\left\vert \kappa_{0}-\kappa_{1}\right\vert }{\left\vert m\right\vert
	^{1/3}}\operatorname*{e}\nolimits^{\kappa_{0}\left(  r_{0}-r_{1}\right)
}\left(  \frac{r_{1}}{r_{0}}\right)  ^{|m|};
\]

\item[d.] for $\left\vert m\right\vert \geq2+r_{1}\max\left\{  \mu\kappa
_{0},\kappa_{1}\right\}  $ it holds%
\begin{equation}
	\left\vert \kappa_{0}A_{m}^{0,2}J_{m}^{\prime}\left(  \kappa_{0}r_{0}\right)
	\right\vert \leq C_{\eta};\label{d_bound}%
\end{equation}

\item[e.] for $\left\vert m\right\vert \geq r_{1}\max\left\{  \mu\kappa
_{0},\kappa_{1}\right\}  $ it holds%
\[
\left\vert A_{m}^{1,2}J_{m}\left(  \kappa_{1}r_{1}\right)  \right\vert \leq
C_{\eta}\frac{\pi r_{0}}{2\left\vert m\right\vert }\left(  1+\frac{2\left\vert
	\kappa_{0}-\kappa_{1}\right\vert r_{1}}{\left\vert m\right\vert ^{1/3}%
}\right)  \operatorname*{e}\nolimits^{\kappa_{0}\left(  r_{0}-r_{1}\right)
}\left(  \frac{r_{1}}{r_{0}}\right)  ^{\left\vert m\right\vert };
\]

\item[f.] for $\left\vert m\right\vert \geq2+r_{1}\max\left\{  \mu\kappa
_{0},\kappa_{1}\right\}  $ it holds%
\[
\left\vert \kappa_{1}A_{m}^{1,2}J_{m}^{\prime}\left(  \kappa_{1}r_{1}\right)
\right\vert \leq C_{\eta}\frac{\pi r_{0}}{2r_{1}}\left(  1+\frac{20\left\vert
	\kappa_{0}-\kappa_{1}\right\vert r_{1}}{3\left\vert m\right\vert ^{1/3}%
}\right)  \operatorname*{e}\nolimits^{\kappa_{0}\left(  r_{0}-r_{1}\right)
}\left(  \frac{r_{1}}{r_{0}}\right)  ^{\left\vert m\right\vert }.
\]

\end{enumerate}
\end{lemma}%

\proof
We prove the bounds only for $m\geq1$ while, for negative $m$, the estimates
follow from (\ref{reflBessel}).

\textbf{@a.}

For $m\geq1$ and $x\in\left[  0,m\right]  $, the estimate%
\begin{equation}
\left\vert J_{m}\left(  x\right)  Y_{m}\left(  x\right)  \right\vert \leq
\frac{2.09}{2\pi\sqrt{m^{2}-x^{2}}}.\label{JYest}%
\end{equation}
follows from \cite[p. 241]{Capdebosq_ball} based on the results in \cite[p.
449]{BoydDunster}. Since the Hankel function can be written in terms of the Bessel
functions of first and second kind:%
\[
H_{m}^{\left(  1\right)  }=J_{m}\left(  x\right)  +\operatorname*{i}%
Y_{m}\left(  x\right),
\]
a triangle inequality leads to%
\[
\left\vert J_{m}\left(  x\right)  H_{m}^{\left(  1\right)  }\left(  x\right)
\right\vert \leq\left\vert J_{m}\left(  x\right)  \right\vert ^{2}+\left\vert
J_{m}\left(  x\right)  Y_{m}\left(  x\right)  \right\vert .
\]
The estimate%

\begin{equation}
|J_{m}\left(  x\right)  |\leq\frac{\left\vert x/2\right\vert ^{m}}{m!}%
\leq\left\vert \frac{\operatorname*{e}x}{2m}\right\vert ^{m}\label{Jmest}%
\end{equation}
for $J_{m}\left(  x\right)  $ follows from \cite[Table 10.14.4]{NIST:DLMF} and
Stirling's formula. In this way, we get%
\begin{equation}
\left\vert J_{m}\left(  x\right)  H_{m}^{\left(  1\right)  }\left(  x\right)
\right\vert \leq\eta^{2m}+\frac{1}{m}\frac{2.09}{2\pi\sqrt{1-\mu^{-2}}}%
\leq\eta^{2m}+\frac{1}{m}\frac{2.09}{2\pi\sqrt{1-\left(  2/\operatorname*{e}%
	\right)  ^{2}}}\leq C_{\eta}m^{-1}\label{JmHmest}%
\end{equation}
and the choice $x=\kappa_{0}r_{0}$ leads to the assertion.

\textbf{@b.}

The coefficient $A_{m}^{0,1}H_{m}^{\left(  1\right)  }\left(  \kappa_{0}%
r_{0}\right)  $ can be written in the form (cf. (\ref{Am01}))%
\[
A_{m}^{0,1}H_{m}^{\left(  1\right)  }\left(  \kappa_{0}r_{0}\right)
=-\frac{\pi\operatorname*{i}r_{0}}{2}S_{m}\left(  \kappa_{0}r_{1}\right)
\frac{J_{m}\left(  \kappa_{0}r_{1}\right)  }{J_{m}\left(  \kappa_{0}%
r_{0}\right)  }J_{m}\left(  \kappa_{0}r_{0}\right)  H_{m}^{\left(  1\right)
}\left(  \kappa_{0}r_{0}\right)  .
\]
From \cite[Prop. 7.1]{Capdebosq_ball} it follows that for $m\geq\max\left\{
\kappa_{0},\kappa_{1}\right\}  r_{1}$ it holds%
\begin{equation}
\left\vert S_{m}\left(  \kappa_{0}r_{1}\right)  \right\vert \leq
2\frac{\left\vert \kappa_{0}-\kappa_{1}\right\vert r_{1}}{m^{1/3}%
}.\label{SMest}%
\end{equation}
The combination with (\ref{JmHmest}) and (\ref{Jmdiffarg}) leads to%
\[
\left\vert A_{m}^{0,1}H_{m}^{\left(  1\right)  }\left(  \kappa_{0}%
r_{0}\right)  \right\vert \leq C_{\eta}\pi r_{0}r_{1}\frac{\left\vert
\kappa_{0}-\kappa_{1}\right\vert }{m^{4/3}}\operatorname*{e}\nolimits^{\kappa
_{0}\left(  r_{0}-r_{1}\right)  }\left(  \frac{r_{1}}{r_{0}}\right)  ^{m}.
\]
\textbf{@c.}

We investigate the Fourier coefficient (cf. (\ref{Am01})):%
\[
\kappa_{0}A_{m}^{0,1}\left(  H_{m}^{\left(  1\right)  }\right)  ^{\prime
}\left(  \kappa_{0}r_{0}\right)  =-\frac{\pi\operatorname*{i}r_{0}\kappa_{0}%
}{2}S_{m}\left(  \kappa_{0}r_{1}\right)  J_{m}\left(  \kappa_{0}r_{1}\right)
\left(  H_{m}^{\left(  1\right)  }\right)  ^{\prime}\left(  \kappa_{0}%
r_{0}\right)  .
\]
In a similar way as in case \textbf{b}, we get%
\begin{equation}
\left\vert \kappa_{0}A_{m}^{0,1}\left(  H_{m}^{\left(  1\right)  }\right)
^{\prime}\left(  \kappa_{0}r_{0}\right)  \right\vert \leq\pi r_{0}r_{1}%
\kappa_{0}\frac{\left\vert \kappa_{0}-\kappa_{1}\right\vert }{m^{1/3}}%
\frac{J_{m}\left(  \kappa_{0}r_{1}\right)  }{J_{m}\left(  \kappa_{0}%
r_{0}\right)  }\left\vert J_{m}\left(  \kappa_{0}r_{0}\right)  \left(
H_{m}^{\left(  1\right)  }\right)  ^{\prime}\left(  \kappa_{0}r_{0}\right)
\right\vert .\label{kappa0c}%
\end{equation}
We use \cite[10.6.1, 10.5.3]{NIST:DLMF} for%
\begin{equation}%
\begin{array}
[c]{ll}%
C_{m}^{\prime}\left(  x\right)  =\frac{1}{2}\left(  C_{m-1}\left(  x\right)
-C_{m+1}\left(  x\right)  \right)   & \text{for }C_{m}\in\left\{  J_{m}%
,Y_{m},H_{m}^{\left(  1\right)  }\right\}  ,\\
J_{m}\left(  x\right)  {H_{m+1}^{(1)}}\left(  x\right)  =J_{m+1}\left(
x\right)  {H_{m}^{(1)}}\left(  x\right)  -2\operatorname*{i}/\left(  \pi
x\right)  , &
\end{array}
\label{crossrelations}%
\end{equation}
and get%
\begin{align}
J_{m}\left(  x\right)  \left(  H_{m}^{\left(  1\right)  }\right)  ^{\prime
}\left(  x\right)   &  =\frac{1}{2}J_{m}\left(  x\right)  \left(
H_{m-1}^{\left(  1\right)  }\left(  x\right)  -H_{m+1}^{\left(  1\right)
}\left(  x\right)  \right)  \label{JmHmprime}\\
&  =\frac{\operatorname*{i}}{\pi x}+\frac{1}{2}\left(  J_{m}\left(  x\right)
H_{m-1}^{\left(  1\right)  }\left(  x\right)  -J_{m+1}\left(  x\right)
{H_{m}^{(1)}}\left(  x\right)  \right)  .\nonumber
\end{align}

Next we compare two consecutive Bessel functions for small argument. We get
from \cite[(4)]{Turan_bessel_baricz} and straightforward monotonicity
properties%
\begin{equation}
\frac{J_{m}\left(  x\right)  }{J_{m-1}\left(  x\right)  }<\frac{m-\sqrt
{m^{2}-\frac{m}{m+1}x^{2}}}{\frac{m}{m+1}x}\leq\frac{m-\sqrt{m^{2}-\frac
	{m}{m+1}\left(  m-1\right)  ^{2}}}{\frac{m}{m+1}\left(  m-1\right)  }%
\leq1\quad\forall x\in\left[  0,m-1\right]  .\label{consecJ}%
\end{equation}
We use (\ref{JmHmprime}) and (\ref{JmHmest}) to get for $m\geq2$ and $\mu
x<m-1:$%
\begin{align*}
\left\vert J_{m}\left(  x\right)  \left(  H_{m}^{\left(  1\right)  }\right)
^{\prime}\left(  x\right)  \right\vert  &  \leq\frac{1}{\pi x}+\frac{1}%
{2}\left(  \left\vert J_{m-1}\left(  x\right)  H_{m-1}^{\left(  1\right)
}\left(  x\right)  \right\vert +\left\vert J_{m}\left(  x\right)  {H_{m}%
^{(1)}}\left(  x\right)  \right\vert \right)  \\
&  \leq\frac{1}{\pi x}+2C_{\eta}m^{-1}.
\end{align*}
The combination with (\ref{kappa0c}) leads to the final estimate (with an
adjusted constant $C_{\eta}$)%
\begin{align*}
\left\vert \kappa_{0}A_{m}^{0,1}\left(  H_{m}^{\left(  1\right)  }\right)
^{\prime}\left(  \kappa_{0}r_{0}\right)  \right\vert  &  \leq\pi r_{0}%
r_{1}\kappa_{0}\frac{\left\vert \kappa_{0}-\kappa_{1}\right\vert }{m^{1/3}%
}\left(  \frac{1}{\pi\kappa_{0}r_{0}}+2\frac{C_{\eta}}{m}\right)  \frac
{J_{m}\left(  \kappa_{0}r_{1}\right)  }{J_{m}\left(  \kappa_{0}r_{0}\right)
}\\
&  \overset{\text{(}m\geq\mu r_{0}\kappa_{0}\text{)}}{\leq}C_{\eta}r_{1}%
\frac{\left\vert \kappa_{0}-\kappa_{1}\right\vert }{m^{1/3}}\frac{J_{m}\left(
\kappa_{0}r_{1}\right)  }{J_{m}\left(  \kappa_{0}r_{0}\right)  }\\
&  \leq C_{\eta}r_{1}\frac{\left\vert \kappa_{0}-\kappa_{1}\right\vert
}{m^{1/3}}\operatorname*{e}\nolimits^{\kappa_{0}\left(  r_{0}-r_{1}\right)
}\left(  \frac{r_{1}}{r_{0}}\right)  ^{m}.
\end{align*}

\textbf{@d. }

The estimate follows in an analogous way as part \textbf{c}:%
\begin{align*}
\left\vert \kappa_{0}A_{m}^{0,2}J_{m}^{\prime}\left(  \kappa_{0}r_{0}\right)
\right\vert  &  \overset{\text{(\ref{Am01})}}{\leq}\frac{\pi\kappa_{0}r_{0}%
}{2}\left\vert H_{m}^{\left(  1\right)  }\left(  \kappa_{0}r_{0}\right)
J_{m}^{\prime}\left(  \kappa_{0}r_{0}\right)  \right\vert \\
&  \overset{\text{(\ref{crossrelations}), (\ref{consecJ})}}{\leq}\frac
{\pi\kappa_{0}r_{0}}{4}\left(  \frac{2}{\pi\kappa_{0}r_{0}}+\left\vert
{H_{m-1}^{(1)}}\left(  \kappa_{0}r_{0}\right)  J_{m-1}\left(  \kappa_{0}%
r_{0}\right)  \right\vert +\left\vert H_{m}^{\left(  1\right)  }\left(
\kappa_{0}r_{0}\right)  J_{m}\left(  \kappa_{0}r_{0}\right)  \right\vert
\right)  \\
&  \leq\frac{1}{2}+\frac{\pi\kappa_{0}r_{0}}{2}C_{\eta}m^{-1}\leq C_{\eta
}^{\operatorname*{III}}.
\end{align*}
By adjusting the constant $C_{\eta}$ the bound (\ref{d_bound}) follows.

\textbf{@e \& f.}

We start with the representations (cf. (\ref{Am01}):%
\begin{align}
A_{m}^{1,2}J_{m}\left(  \kappa_{1}r_{1}\right)   &  =\frac{\pi
\operatorname*{i}r_{0}}{2}\left(  1-S_{m}\left(  \kappa_{0}r_{1}\right)
\right)  J_{m}\left(  \kappa_{0}r_{1}\right)  H_{m}^{\left(  1\right)
}\left(  \kappa_{0}r_{0}\right)  ,\label{fouriecoeffE}\\
\kappa_{1}A_{m}^{1,2}J_{m}^{\prime}\left(  \kappa_{1}r_{1}\right)   &
=\frac{\pi\operatorname*{i}r_{0}}{2}\kappa_{0}\left(  \frac{J_{m}^{\prime
}\left(  \kappa_{0}r_{1}\right)  }{J_{m}\left(  \kappa_{0}r_{1}\right)
}-S_{m}\left(  \kappa_{0}r_{1}\right)  \frac{\left(  H_{m}^{\left(  1\right)
}\right)  ^{\prime}\left(  \kappa_{0}r_{1}\right)  }{H_{m}^{\left(  1\right)
}\left(  \kappa_{0}r_{1}\right)  }\right)  J_{m}\left(  \kappa_{0}%
r_{1}\right)  H_{m}^{\left(  1\right)  }\left(  \kappa_{0}r_{0}\right)
.\label{ratiobesselhankel}%
\end{align}
To estimate the first term we employ (\ref{SMest}), (\ref{Jmdiffarg}),
(\ref{JmHmest}) and get%
\[
\left\vert A_{m}^{1,2}J_{m}\left(  \kappa_{1}r_{1}\right)  \right\vert \leq
C_{\eta}\frac{\pi r_{0}}{2m}\left(  1+\frac{2\left\vert \kappa_{0}-\kappa
_{1}\right\vert r_{1}}{m^{1/3}}\right)  \operatorname*{e}\nolimits^{\kappa
_{0}\left(  r_{0}-r_{1}\right)  }\left(  \frac{r_{1}}{r_{0}}\right)  ^{m}.
\]
For the second coefficient (\ref{ratiobesselhankel}) we employ (\ref{SMest}),
(\ref{Propratios})%
\begin{align*}
\left\vert \kappa_{1}A_{m}^{1,2}J_{m}^{\prime}\left(  \kappa_{1}r_{1}\right)
\right\vert  & \leq\frac{\pi\kappa_{0}r_{0}}{2}\left(  1+\frac{20}{3}%
\frac{\left\vert \kappa_{0}-\kappa_{1}\right\vert r_{1}}{m^{1/3}}\right)
\left|\frac{J_{m}^{\prime}\left(  \kappa_{0}r_{1}\right)  }{J_{m}\left(  \kappa
_{0}r_{1}\right)  }\frac{J_{m}\left(  \kappa_{0}r_{1}\right)  }{J_{m}\left(
\kappa_{0}r_{0}\right)  }\right|\left|J_{m}\left(  \kappa_{0}r_{0}\right)  H_{m}^{\left(
1\right)  }\left(  \kappa_{0}r_{0}\right)\right|  \\
& \leq C_{\eta}\frac{\pi\kappa_{0}r_{0}}{2}\left(  1+\frac{20}{3m^{1/3}%
}\left\vert \kappa_{0}-\kappa_{1}\right\vert r_{1}\right)  \frac{1}{\kappa
_{0}r_{1}}\operatorname*{e}\nolimits^{\kappa_{0}\left(  r_{0}-r_{1}\right)
}\left(  \frac{r_{1}}{r_{0}}\right)  ^{m}\\
& \leq C_{\eta}\frac{\pi r_{0}}{2r_{1}}\left(  1+\frac{20\left\vert \kappa
_{0}-\kappa_{1}\right\vert r_{1}}{3m^{1/3}}\right)  \operatorname*{e}%
\nolimits^{\kappa_{0}\left(  r_{0}-r_{1}\right)  }\left(  \frac{r_{1}}{r_{0}%
}\right)  ^{m}.
\end{align*}%
\endproof

\subsubsection{Spectral approximability of the solution}

For a function $f\in H^{s}\left(  \Gamma_{j}\right)  $, $s\in\mathbb{R}$, we
denote its complex Fourier coefficients by%
\[
f_{m}:=\frac{1}{r_{j}}{\color{black}\langle  f,\operatorname*{e}\nolimits^{\operatorname*{i}%
m\theta}\rangle  _{\Gamma_{j}}}=\int_{-\pi}^{\pi}f(
r\operatorname*{e}\nolimits^{\operatorname*{i}\theta})
\operatorname*{e}\nolimits^{-\operatorname*{i}m\theta}d\theta
\]
so that%
\[
f(r\operatorname*{e}\nolimits^{\operatorname*{i}\theta})
=\frac{1}{2\pi}\sum_{m\in\mathbb{Z}}f_{m}g_{m},
\]
{\color{black}with $g_m$ as in \eqref{Defgm}}.
We define its norm by%
\[
\left\Vert f\right\Vert _{H^{s}\left(  \Gamma_{j}\right)  }:=\left(
\frac{r_{j}}{2\pi}\sum_{m\in\mathbb{Z}}\left(  1+m^{2}\right)  ^{s}\left\vert
f_{m}\right\vert ^{2}\right)  ^{1/2}.
\]
For given $M\in\mathbb{N}$, we will estimate the approximation properties of
the spectral space%
\[
S_{M}^{j}:=\operatorname*{span}\left\{  \left.  g_{m}\right\vert _{\Gamma_{j}%
}:m\in\mathbb{Z}\wedge\left\vert m\right\vert \leq M\right\}  .
\]
To this aim, let 
\[%
\begin{array}
[c]{cc}%
e_{\operatorname*{D};0}^{0}\left(  M\right)  :=\inf_{v\in S_{M}^{0}}\left\Vert
u_{\operatorname*{D};0}^{0}-v\right\Vert _{H^{1/2}\left(  \Gamma_{0}\right)
}, & e_{\operatorname*{N};0}^{0}\left(  M\right)  :=\inf_{v\in S_{M}^{0}%
}\left\Vert u_{\operatorname*{N};0}^{0}-v\right\Vert _{H^{-1/2}\left(
\Gamma_{0}\right)  },\\
e_{\operatorname*{D};1}\left(  M\right)  :=\inf_{v\in S_{M}^{1}}\left\Vert
u_{\operatorname*{D};1}-v\right\Vert _{H^{1/2}\left(  \Gamma_{1}\right)  }, &
e_{\operatorname*{N};1}\left(  M\right)  :=\inf_{v\in S_{M}^{1}}\left\Vert
u_{\operatorname*{N};1}-v\right\Vert _{H^{-1/2}\left(  \Gamma_{1}\right)  }.
\end{array}
\]

The Sobolev regularity of the solution $u$ of (\ref{transmission_problem_1})--(\ref{transmission_problem_4}) are well
studied. For given $g\in H^{s}\left(  \Gamma_{0}\right)  $ for some
$s\geq-1/2$ the traces of the solution on $\Gamma_{0}$ satisfy%
\[
u_{\operatorname*{D};0}^{0}\in H^{s+1}\left(  \Gamma_{0}\right)  ,\quad
u_{\operatorname*{N};0}^{0}\in H^{s}\left(  \Gamma_{0}\right)  ,
\]
while the other traces $u_{\operatorname*{D};1}^{1}$ and $u_{\operatorname*{N}%
;0}^{1}$ are analytic.

\begin{theorem}\label{th:estimates}
For $\eta\in\left]  0,1\right[  $, let $C_{\eta}$ be as in Lemma \ref{LemFC}
and $\mu>\operatorname*{e}/\left(  2\eta\right)  $. Let $M\geq2+r_{0}\mu\max\left\{  \kappa_{0},\kappa_{1}\right\}  $ 
and $\nu:=\frac{1}{2}\left(
\frac{r_{0}-r_{1}}{r_{0}}\right)  \left(  \frac{\mu-1}{\mu}\right)  >0$. For
given $g\in H^{s}\left(  \Gamma_{0}\right)  $ for some $s\geq-1/2$, it holds%
\begin{align*}
e_{\operatorname*{D};0}^{0}\left(  M\right)   &  \leq C_{\eta}\sqrt{2}%
r_{0}\left(  \frac{1}{\nu}+1\right)  \left(  1+M\right)  ^{-1/2-s}\left\Vert
{\color{black}g}\right\Vert _{H^{s}\left(  \Gamma_{0}\right)  },\\
e_{\operatorname*{N};0}^{0}\left(  M\right)   &  \leq C_{\eta}\left(  \frac
{1}{\nu}+1\right)  \left(  1+M\right)  ^{-1/2-s}\left\Vert g\right\Vert
_{H^{s}\left(  \Gamma_{0}\right)  },\\
e_{\operatorname*{D};1}\left(  M\right)   &  \leq C_{\eta}\pi\sqrt{\frac
	{r_{0}r_{1}}{2}}\left(  \frac{2}{\nu}+1\right)  \operatorname*{e}%
\nolimits^{-\nu {\color{black}M}}\left(  1+M\right)  ^{-1/2-s}\left\Vert g\right\Vert
_{H^{s}\left(  \Gamma_{0}\right)  },\\
e_{\operatorname*{N};1}\left(  M\right)   &  \leq C_{\eta}\pi\sqrt{\frac
	{r_{0}}{2r_{1}}}\left(  \frac{20}{3\nu}+1\right)  \operatorname*{e}%
\nolimits^{-\nu {\color{black}M}}\left(  1+{\color{black}M}\right)  ^{-1/2-s}\left\Vert g\right\Vert
_{H^{s}\left(  \Gamma_{0}\right)  }.
\end{align*}

\end{theorem}%
\proof
For the error $e_{\operatorname*{D};0}^{0}\left(  M\right)  $, by employing @a and @b of Lemma \ref{lm:estimates_four_coeff}, we obtain%
\begin{align*}
& e_{\operatorname*{D};0}^{0}\left(  M\right)  \leq\sqrt{\frac{r_{0}}{2\pi}%
}\left(  \sum_{\left\vert m\right\vert > M}\left(  1+m^{2}\right)
^{1/2-s}\left(  \left\vert A_{m}^{0,1}H_{m}^{\left(  1\right)  }\left(
\kappa_{0}r_{0}\right)  \right\vert +\left\vert A_{m}^{0,2}J_{m}\left(
\kappa_{0}r_{0}\right)  \right\vert \right)  ^{2}\left(  1+m^{2}\right)
^{s}\left\vert \alpha_{m}\right\vert ^{2}\right)  ^{1/2}\\
& \quad\leq C_{\eta}\sqrt{\frac{r_{0}}{2\pi}}\left(  \sum_{\left\vert
m\right\vert > M}\frac{\left(  1+m^{2}\right)  ^{1/2-s}}{m^{2}}\left(
r_{0}r_{1}\frac{\left\vert \kappa_{0}-\kappa_{1}\right\vert }{{\color{black}|m|^{1/3}}%
}\operatorname*{e}\nolimits^{\kappa_{0}\left(  r_{0}-r_{1}\right)  }\left(
\frac{r_{1}}{r_{0}}\right)  ^{{\color{black}|m|}}+r_{0}\right)  ^{2}\left(  1+m^{2}\right)
^{s}\left\vert \alpha_{m}\right\vert ^{2}\right)  ^{1/2}\\
& \quad\overset{\text{(\ref{combineexpest})}}{\leq}C_{\eta}\sqrt{\frac{r_{0}%
}{2\pi}}r_{0}\left(  \sum_{\left\vert m\right\vert > M}\frac{\left(
1+m^{2}\right)  ^{1/2-s}}{m^{2}}\left(  \nu^{-1}\operatorname*{e}%
\nolimits^{-\nu {\color{black}|m|}}+1\right)  ^{2}\left(  1+m^{2}\right)  ^{s}\left\vert
\alpha_{m}\right\vert ^{2}\right)  ^{1/2}\\
& \quad\leq C_{\eta}\sqrt{2}r_{0}\left(  \nu^{-1}+1\right)  \left(
1+M\right)  ^{-1/2-s}\left\Vert g\right\Vert _{H^{s}\left(  \Gamma
_{0}\right)  }.
\end{align*}
with $\nu$ as in the statement of Lemma \ref{lm:lemma_A1}.

For the error $e_{\operatorname*{N};0}^{0}\left(  M\right)  $, it holds%
\begin{align*}
e_{\operatorname*{N};0}^{0}\left(  M\right)    & \leq C_{\eta}\sqrt
{\frac{r_{0}}{2\pi}}\left(  \sum_{\left\vert m\right\vert > M}\left(
1+m^{2}\right)  ^{-1/2}\left(  \frac{\left\vert \kappa_{0}-\kappa
_{1}\right\vert r_{1}}{{\color{black}\left\vert m\right\vert ^{1/3}}}\operatorname*{e}%
\nolimits^{\kappa_{0}\left(  r_{0}-r_{1}\right)  }\left(  \frac{r_{1}}{r_{0}%
}\right)  ^{{\color{black}\left\vert m\right\vert} }+1\right)  ^{2}\left\vert \alpha
_{m}\right\vert ^{2}\right)  ^{1/2}\\
& \leq C_{\eta}\sqrt{\frac{r_{0}}{2\pi}}\left(  \sum_{\left\vert m\right\vert
> M}\left(  1+m^{2}\right)  ^{-1/2}\left(  \nu^{-1}\operatorname*{e}%
\nolimits^{-\nu {\color{black}|m|}}+1\right)  ^{2}\left\vert \alpha_{m}\right\vert ^{2}\right)
^{1/2}.%
\end{align*}
Reasoning as in the previous case leads to%
\[
e_{\operatorname*{N};0}^{0}\left(  M\right)  \leq C_{\eta}\left(  \nu
^{-1}+1\right)  \left(  1+M\right)  ^{-1/2-s}\left\Vert g\right\Vert
_{H^{s}\left(  \Gamma_{0}\right)  }.
\]

Next we estimate the error $e_{\operatorname*{D};1}\left(  M\right)  :$%
\begin{align*}
& e_{\operatorname*{D};1}\left(  M\right)  \leq\sqrt{\frac{r_{1}}{2\pi}%
}\left(  \sum_{\left\vert m\right\vert > M}\left(  1+m^{2}\right)
^{1/2-s}\left\vert A_{m}^{1,2}J_{m}\left(  \kappa_{1}r_{1}\right)  \right\vert
^{2}\left(  1+m^{2}\right)  ^{s}\left\vert \alpha_{m}\right\vert ^{2}\right)
^{1/2}\\
& \quad\leq C_{\eta}\frac{\pi r_{0}}{2}\sqrt{\frac{r_{1}}{2\pi}}\times\\
& \quad\times\left(  \sum_{\left\vert m\right\vert > M}\frac{\left(
1+m^{2}\right)  ^{1/2-s}}{m^{2}}\left(  \left(  1+\frac{2\left\vert \kappa
_{0}-\kappa_{1}\right\vert r_{1}}{{\color{black}\left\vert m\right\vert ^{1/3}}}\right)
\operatorname*{e}\nolimits^{\kappa_{0}\left(  r_{0}-r_{1}\right)  }\left(
\frac{r_{1}}{r_{0}}\right)  ^{{\color{black}\left\vert m\right\vert} }\right)  ^{2}\left(
1+m^{2}\right)  ^{s}\left\vert \alpha_{m}\right\vert ^{2}\right)  ^{1/2}.
\end{align*}
The combination of this with estimates (\ref{exp2nue}) and
(\ref{combineexpest}) leads to%
\begin{align*}
e_{\operatorname*{D};1}\left(  M\right)    & \leq C_{\eta}\frac{\pi r_{0}}%
{2}\sqrt{\frac{r_{1}}{2\pi}}\left(  \sum_{\left\vert m\right\vert > M}%
\frac{\left(  1+m^{2}\right)  ^{-1/2-s}}{m^{2}}\left(  \left(
\operatorname*{e}\nolimits^{-2\nu {\color{black}|m|}}+\frac{2}{\nu}\operatorname*{e}%
\nolimits^{-\nu {\color{black}|m|}}\right)  \right)  ^{2}\left(  1+m^{2}\right)  ^{s}\left\vert
\alpha_{m}\right\vert ^{2}\right)  ^{1/2}\\
& \leq C_{\eta}\pi\sqrt{\frac{r_{0}r_{1}}{2}}\left(  \frac{2}{\nu}+1\right)
\operatorname*{e}\nolimits^{-\nu M}\left(  1+M\right)  ^{-1/2-s}\left\Vert
g\right\Vert _{H^{s}\left(  \Gamma_{0}\right)  }.
\end{align*}

Finally, for the approximation of the Neumann trace on $\Gamma_{1}$ we get%
\begin{align*}
& e_{\operatorname*{N};1}\left(  M\right)  \leq\sqrt{\frac{r_{1}}{2\pi}%
}\left(  \sum_{\left\vert m\right\vert > M}\left(  1+m^{2}\right)
^{-1/2}\left\vert \kappa_{1}A_{m}^{1,2}J_{m}^{\prime}\left(  \kappa_{1}%
r_{1}\right)  \right\vert ^{2}\left\vert \alpha_{m}\right\vert ^{2}\right)
^{1/2}\\
& \quad\leq\sqrt{\frac{r_{1}}{2\pi}}C_{\eta}\frac{\pi r_{0}}{2r_{1}}\times\\
& \quad\quad\times\left(  \sum_{\left\vert m\right\vert > M}\left(
1+m^{2}\right)  ^{-1/2}\left\vert \left(  1+\frac{20\left\vert \kappa
_{0}-\kappa_{1}\right\vert r_{1}}{3{\color{black}\left\vert m\right\vert ^{1/3}}}\right)
\operatorname*{e}\nolimits^{\kappa_{0}\left(  r_{0}-r_{1}\right)  }\left(
\frac{r_{1}}{r_{0}}\right)  ^{{\color{black}\left\vert m\right\vert} }\right\vert
^{2}\left\vert \alpha_{m}\right\vert ^{2}\right)  ^{1/2}.
\end{align*}
Similarly as for $e_{\operatorname*{D};1}\left(  M\right)  $ we conclude that%
\[
e_{\operatorname*{N};1}\left(  M\right)  \leq C_{\eta}\pi\sqrt{\frac{r_{0}%
}{2r_{1}}}\left(  \frac{20}{3\nu}+1\right)  \operatorname*{e}\nolimits^{-\nu
{\color{black}M}}\left(  1+{\color{black}M}\right)  ^{-1/2-s}\left\Vert g\right\Vert _{H^{s}\left(
\Gamma_{0}\right)  }.
\]%
\endproof

\section{Numerical results}\label{sec:numerical_results}
{\color{black}The aim of this section is twofold: first
we demonstrate the behaviour of the solutions \eqref{useriea}--\eqref{useried} in the case of frequencies close to resonances, and then we validate the convergence estimates provided in Theorem \ref{th:estimates}}.
To this aim, we fix the values of the parameters $r_1 = 0.5$, $n_0 = 0.5$ and $n_1 = 1$, and we consider values of $m$ and $\kappa$ close to the zero of \eqref{eq:wronskian}.
To retrieve the numerical solution, we solve the linear system 

\newcommand{\uDooM}{u_{D;0}^{0,M}}
\newcommand{\uNooM}{u_{N;0}^{0,M}}
\newcommand{\uDjM}{u_{D;1}^M}
\newcommand{\uNjM}{u_{N;1}^M}
\newcommand{\uDiM}{u_{D;1}^M}
\newcommand{\uNiM}{u_{N;1}^M}

\begin{eqnarray}\label{eq:lin_syst}
	&	\hskip-2cm\left[
	\begin{array}{cccc}
		-\WWounouno - \WWuno & \KKpuno -\KKpounouno  & -\WWounoo & \KKpounoo \\
		\\
		-\KKounouno+\KKuno & -\VVounouno-\VVuno  & -\KKounoo & \VVounoo  \\
		\\
		-\WWoouno & -\KKpoouno & -\WWooo & \frac{1}{2}{\mathbb D}+\KKpooo \\
		\\
		\KKoouno & \VVoouno & 2\KKooo & -2\VVooo \\
	\end{array}
	\right]
	\, \left[
	\begin{array}{c}
		\uuDuno\\
		\\
		\uuNuno\\
		\\
		\uuDoo\\
		\\
		\uuNoo\\		
	\end{array}
	\right] 
	&=
	\left[
	\begin{array}{c}
		\bo\\
		\\
		\bo\\
		\\
		\bo\\
		\\
		\VVooo \bg^M\\
	\end{array}
	\right],
\end{eqnarray}
{\color{black}where $\mathbb{M\in}\left\{  \mathbb{M}_{k}^{j,\ell},\mathbb{M}_{1}%
\mid\mathbb{M}\in\left\{  \mathbb{V},\mathbb{K},\mathbb{K}^{\prime}%
,\mathbb{W}\right\}  ,\left( k,j, \ell\right)  \in\mathcal{J}\right\}  $ with
$\mathcal{J}:=\left\{  \left(  0,0,0\right)  ,\left(  0,1,0\right)  ,\left(
0,0,1\right)  ,\left(  0,1,1\right)  \right\}  $ denote the representation of
the Galerkin discretization with respect to the spectral basis $\left\{
g_{m}:\left\vert m\right\vert \leq M\right\}  $ of the corresponding operator
$M\in\left\{  M_{k}^{j,\ell},M_{1}\mid M\in\left\{  V,K,K^{\prime},W\right\}
,\left(k,j,  \ell\right)  \in\mathcal{J}\right\}  $ by diagonal matrices (see
Lemma \ref{lm:potentials_diag})}:
$$[\mathbb{M}_1]_{m,m} = \langle g_m,M_1 g_m\rangle_{\Gamma_1} \qquad |m| \leq M$$
and, for $j,\ell = 0,1$,
$$[\MMo^{j,\ell}]_{m,m} = \langle g_m,M_0^{j,\ell} g_m\rangle_{\Gamma_j} \qquad |m| \leq M,$$
and the right-hand side vector by  $\bg^M = \left[\alpha_m\right]_{|m|\leq M}$ (see \eqref{grep}).
It follows that the solution vectors 
\begin{equation*}
\begin{aligned}
&\displaystyle{\uuDuno = \left[({\uDiM})_m\right]_{|m|\leq M}}, &\displaystyle{\uuNuno = \left[({\uNiM})_m\right]_{|m|\leq M}},\\
&&\\
&\uuDoo =\left[({\uDooM})_m\right]_{|m|\leq M}, &\uuNoo = \left[({\uNooM})_m\right]_{|m|\leq M}
\end{aligned}
\end{equation*}
are such that the approximations of \eqref{useriea}--\eqref{useried} are given by
$$\uDuno(\x) \approx \uDiM(\x) = \sum_{|m|\leq M} ({\uDiM})_m g_m(\x), \qquad \uNuno(\x) \approx \uNiM(\x) = \sum_{|m|\leq M} ({\uNiM})_m g_m(\x),\quad \x\in\Gamma_1,$$
$$\uDoo(\x) \approx \uDooM(\x) = \sum_{|m|\leq M} ({\uDooM})_m g_m(\x), \qquad \uNoo(\x) \approx \uNooM(\x) = \sum_{|m|\leq M} ({\uNooM})_m g_m(\x),\quad \x\in\Gamma_0.$$
To retrieve the solution of Problem \eqref{transmission_problem_1}--\eqref{transmission_problem_4} in the whole domain $\Omega = \Omega_1\cup\Omega_0$, we use the representation formulas
\begin{equation*}
u_1(\x) = -\int_{\Gamma_1} G_{\kappa_1}(\x,\y) \uNuno(\y) \,\text{d}\y + \int_{\Gamma_1} \left\langle  \mathbf{n}_{1}\left(  \mathbf{y}\right)
,\nabla_{\mathbf{y}}G_{\kappa_{1}}\left(  \mathbf{x},\mathbf{y}\right)
\right\rangle \uDuno(\y)\,\text{d}\y , \qquad \x\in \Omega_1,
\end{equation*}
and
\begin{eqnarray*}
u_0(\x) &=& \int_{\Gamma_1} G_{\kappa_0}(\x,\y) \uNuno(\y) \,\text{d}\y + \int_{\Gamma_1} \left\langle  \mathbf{n}_{0}\left(  \mathbf{y}\right)
,\nabla_{\mathbf{y}}G_{\kappa_{0}}\left(  \mathbf{x},\mathbf{y}\right)
\right\rangle \uDuno(\y)\,\text{d}\y \\
 &&-\int_{\Gamma_0} G_{\kappa_0}(\x,\y) \uNoo(\y) \,\text{d}\y + \int_{\Gamma_0} \left\langle  \mathbf{n}_{0}\left(  \mathbf{y}\right)
 ,\nabla_{\mathbf{y}}G_{\kappa_{0}}\left(  \mathbf{x},\mathbf{y}\right)
 \right\rangle \uDoo(\y)\,\text{d}\y  , \qquad \x\in \Omega_0.
\end{eqnarray*}

\newcommand{\bm}{{\overline{m}}}
\paragraph{Example 1.}
% il dato per il test è g2 con m=40
 We consider the datum $g(\x) = g_\bm(\x)$ (see \eqref{Defgm}), with $\bm = 40$. For this choice, the resonance $\kappa^*$, the root of \eqref{eq:wronskian}, is approximately given by the complex number $\kappa^* \approx 90.11-\operatorname*{i}\, 9.64\cdot 10^{-4}$. 
 We solve the problem for a frequency value close to the resonance, specifically $\kappa = \mathfrak{R}(\kappa^*)$, and compute the solution of the linear system \eqref{eq:lin_syst} by selecting $M = \overline{m}$. With this choice, as expected, we achieve full accuracy in approximating the four traces.
 
\begin{figure}[h!]
	\begin{center}
		\includegraphics[height=3.5cm,width=9cm]{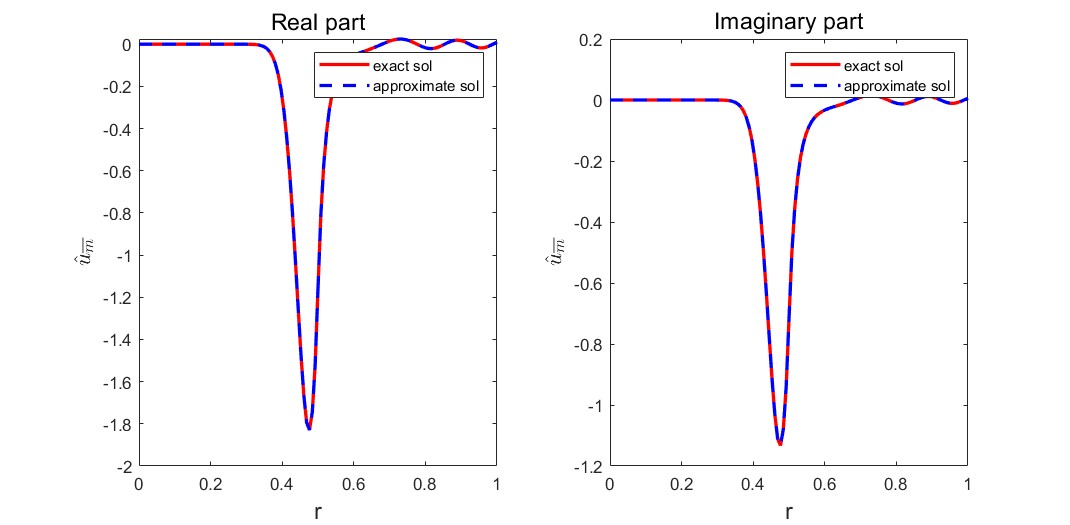}
		\includegraphics[height=3.5cm,width=9cm]{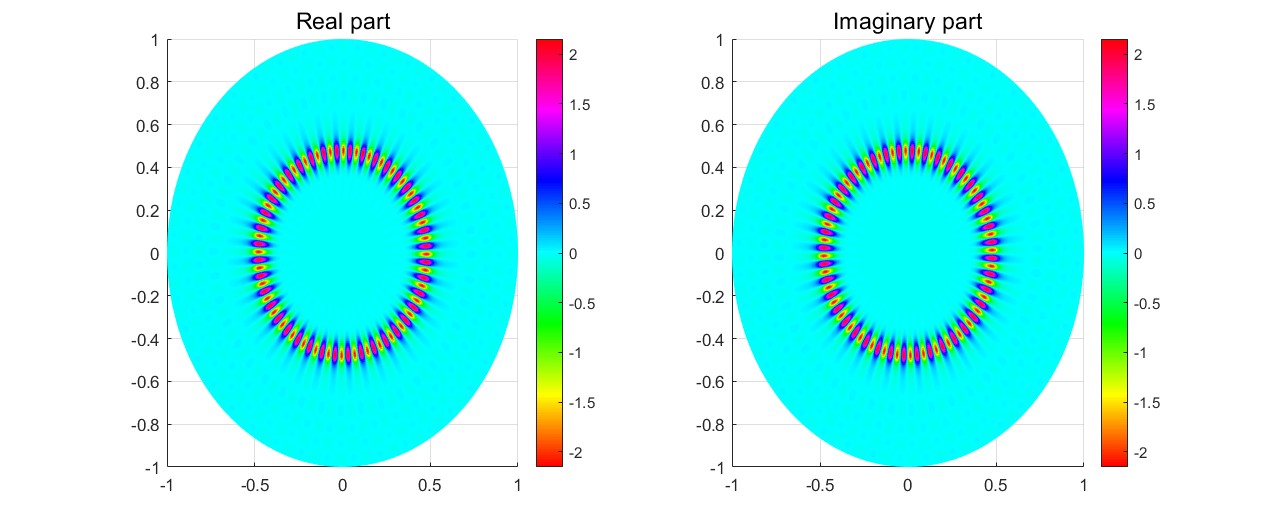}	
		\caption{\small Example 1. Behaviour of the real and imaginary part of the solution in the radial direction (first two left plots) and in the whole domain (last two right plots), for $\bm = 40$ and $\kappa \approx 90.11$.}
		\label{fig:ex_1_Ref_sol_m40_radial}
	\end{center}
\end{figure}

 The left panel of Figure \ref{fig:ex_1_Ref_sol_m40_radial} compares the exact and approximate radial solutions $\hat{u}_{\bm}$ for $r\in[0,1]$. In the right panel, we display the solution over the entire domain $\Omega$. As shown, the solution is localized near $\Gamma_1$, indicating that frequencies close to resonance produce a field concentrated in the vicinity of $\Gamma_1$. {\color{black}This interference phenomenon is commonly referred to in the literature as the whispering gallery mode.}

\begin{figure}[h!]
	\begin{center}
		\includegraphics[height=3.5cm,width=9cm]{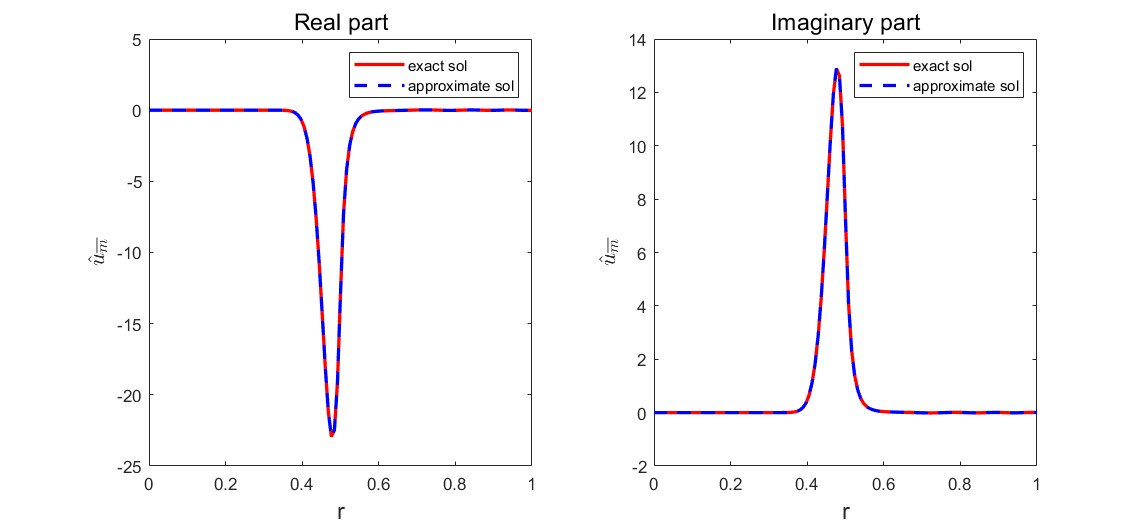}
		\includegraphics[height=3.5cm,width=9cm]{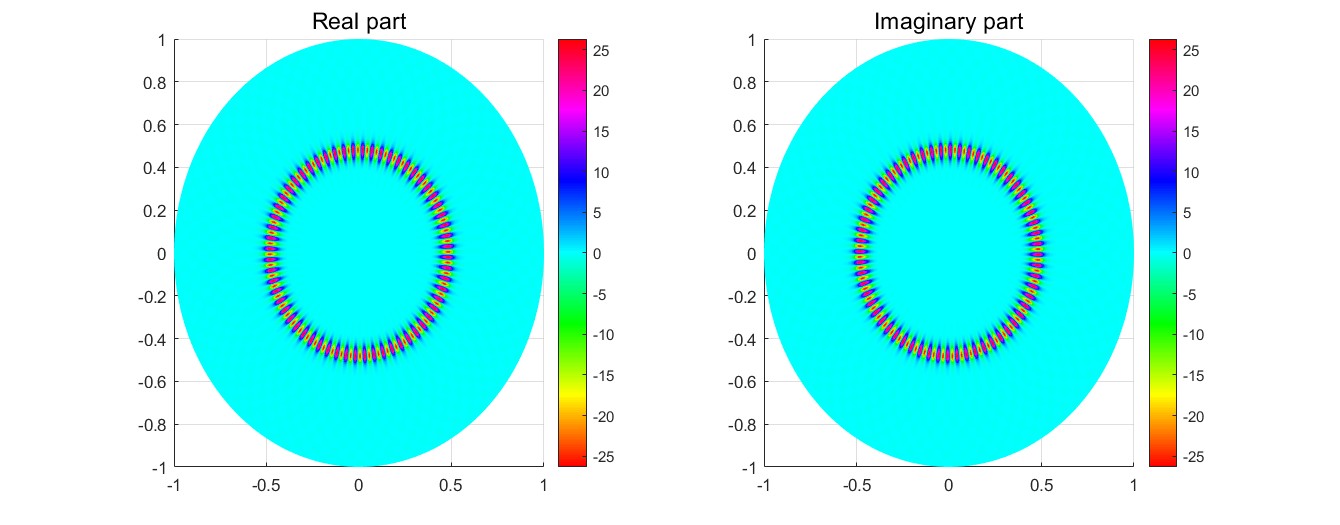}
		\caption{\small Example 2. Behaviour of the real and imaginary part of the solution in the radial direction (first two left plots) and in the whole domain (last two right plots), for $\bm = 60$ and $\kappa \approx 131.97$.}
		\label{fig:ex_2_Ref_sol_m60_radial}
	\end{center}
\end{figure}

\paragraph{Example 2.}
% il dato per il test è g2 con m=60
For the higher mode value $\bm = 60$, the resonance frequency is approximately $\kappa^* \approx 131.97 - \operatorname*{i} \, 3.87\cdot 10^{-6}$.
Once again, we solve the problem for the frequency $\kappa = \mathfrak{R}(\kappa^*)$, obtaining the approximate solution by setting $M = \overline{m}$. As shown in Figure \ref{fig:ex_2_Ref_sol_m60_radial}, the localization of the solution becomes increasingly pronounced near the interface, with its maximum value growing, as the frequency $\kappa$ increases.

\clearpage
\paragraph{Example 3.}
% il dato per il test è g6
In this test we consider the value of the frequency $\kappa \approx  90.11$, which is close to the resonance $\kappa^*$ of Example 1. We choose the datum $g(\x) = |x_2| \in H^{s}(\Gamma_0)$, with $s = 3/2-\varepsilon$ for any $\varepsilon > 0$. 
To validate the results of Theorem \ref{th:estimates}, being the exact solution unknown, we compute the \emph{reference} traces $\uDuno, \uNuno, \uDo$ and $\uNo$ on $\Gamma_1$ and $\Gamma_0$, respectively, associated with the choice of the reference value $\overline{M} = 300$. 

\begin{figure}[h!]
	\begin{center}
		\includegraphics[height=4cm,width=8cm]{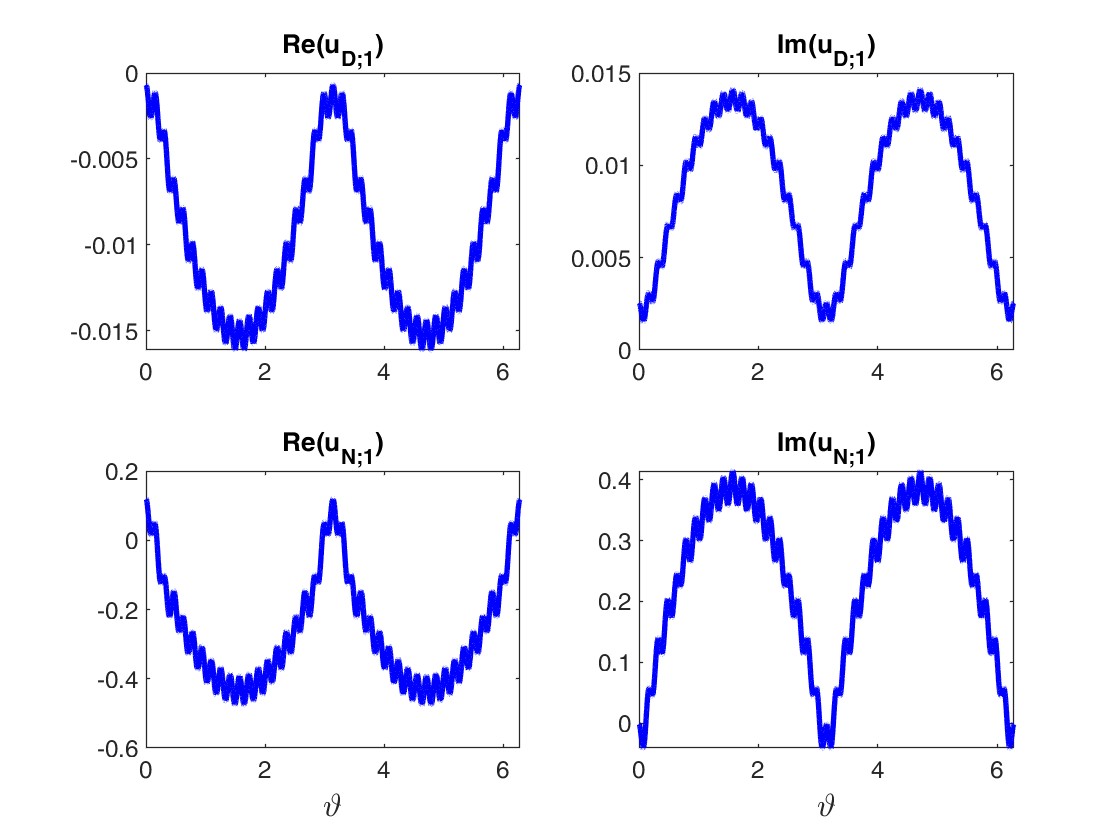}\includegraphics[height=4cm,width=8cm]{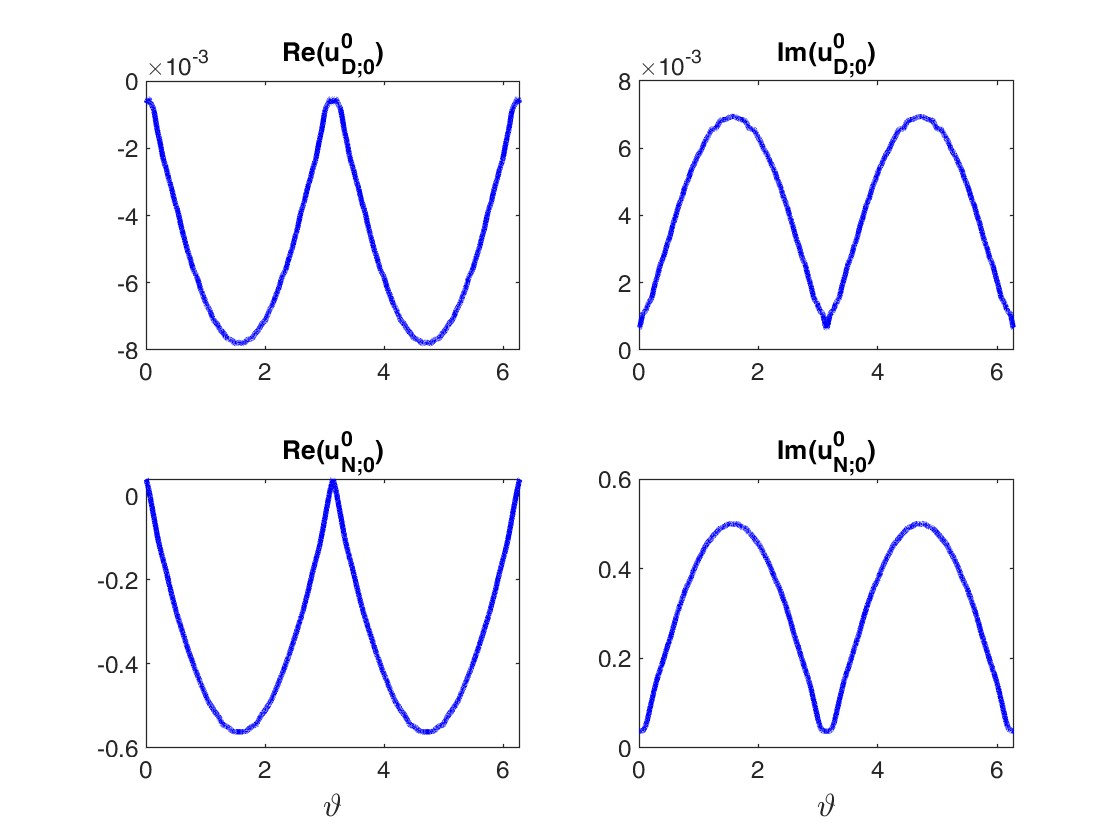}	
		
		\includegraphics[height=3cm,width=8cm]{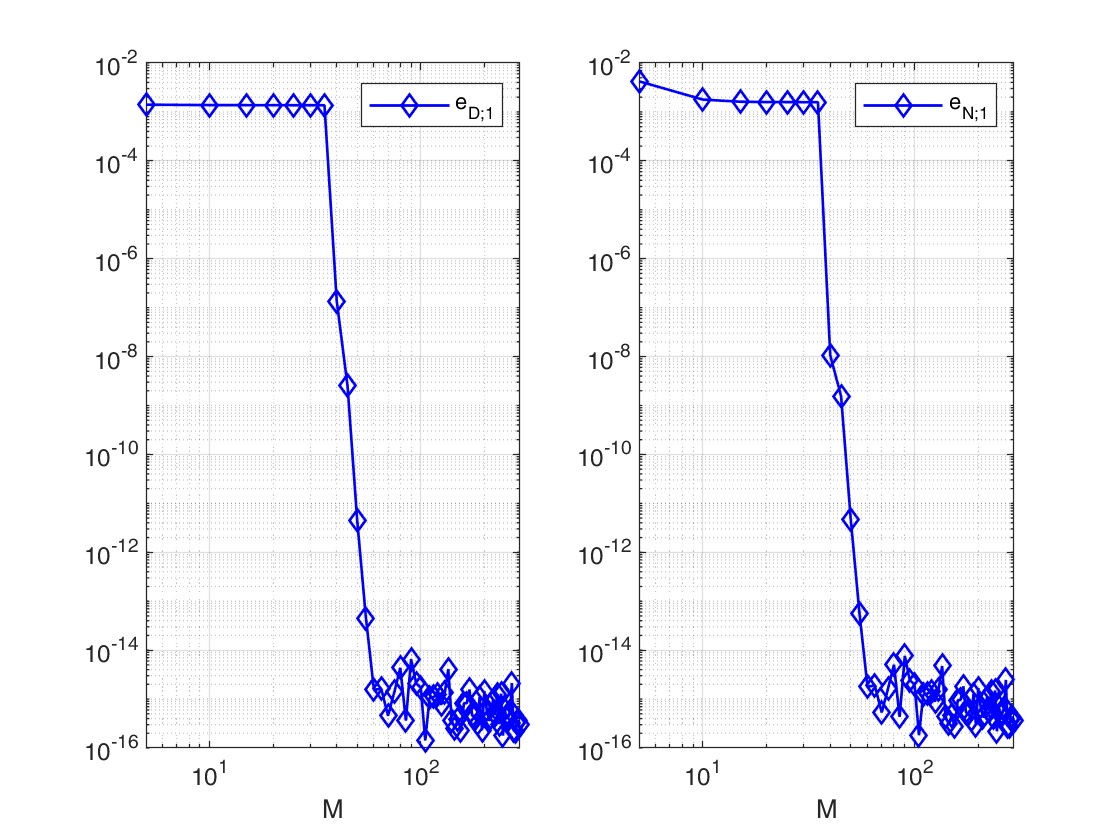}\includegraphics[height=3cm,width=8cm]{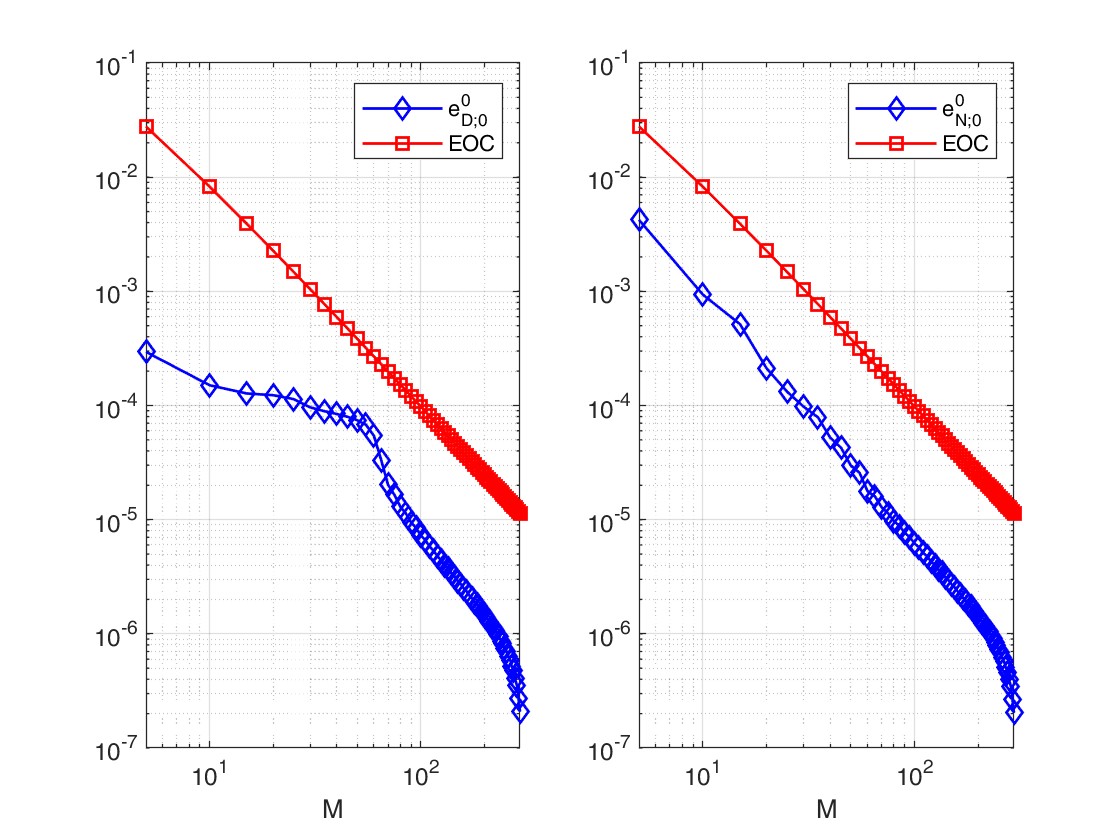}			
		\caption{Example 3. Behaviour of the traces $\uDuno, \uNuno$ (first two rows and columns) and $\uDo, \uNo$ (first two rows and second and last columns) with respect to the angle $\vartheta$. Log-log scale of the trace errors (last row).}
		\label{fig:ex3_traces}
	\end{center}
\end{figure}

In Figure \ref{fig:ex3_traces} we represent the real and imaginary parts of the four reference traces behaviour, with respect to the angle parameter $\vartheta\in[0,2\pi]$. We also report the behaviour of the $H^{1/2}$-norm and $H^{-1/2}$-norm of the Dirichlet and Neumann trace errors, respectively, with respect to increasing values of the discretization parameter $M$. The numerical results confirm the error estimates provided by Theorem \ref{th:estimates}. Indeed, for what concerns the errors $e_{\operatorname*{D};1}$ and $e_{\operatorname*{N};1}$, for $M$ large enough, we observe an exponential decay, till the machine precision is reached. For what concerns $e_{\operatorname*{D};0}$ and $e_{\operatorname*{N};0}$, both errors scale as the expected ones, as the comparison with the graph of the EOC confirms.

\paragraph{Example 4.}
% il dato per il test è g5
In the same frequency regime of Example 3, we consider now the datum $g(\x) = g_\bm(\x) + e^{x_1} \in H^{s}(\Gamma_0)$, for any $s > 0$.
The \emph{reference} traces $\uDuno, \uNuno, \uDo$ and $\uNo$ on $\Gamma_1$ and $\Gamma_0$, have been computed by choosing $\overline{M} = 300$. 

\begin{figure}[h!]
	\begin{center}
			\includegraphics[height=4cm,width=8cm]{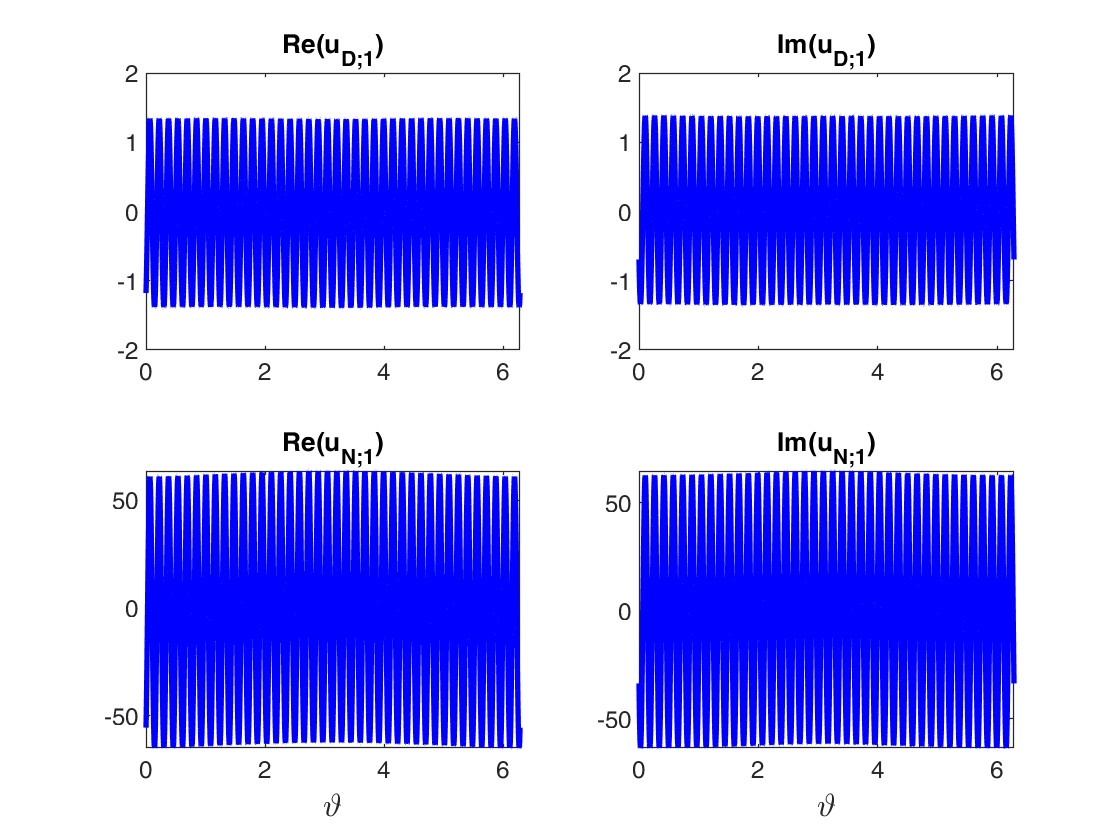}\includegraphics[height=4cm,width=8cm]{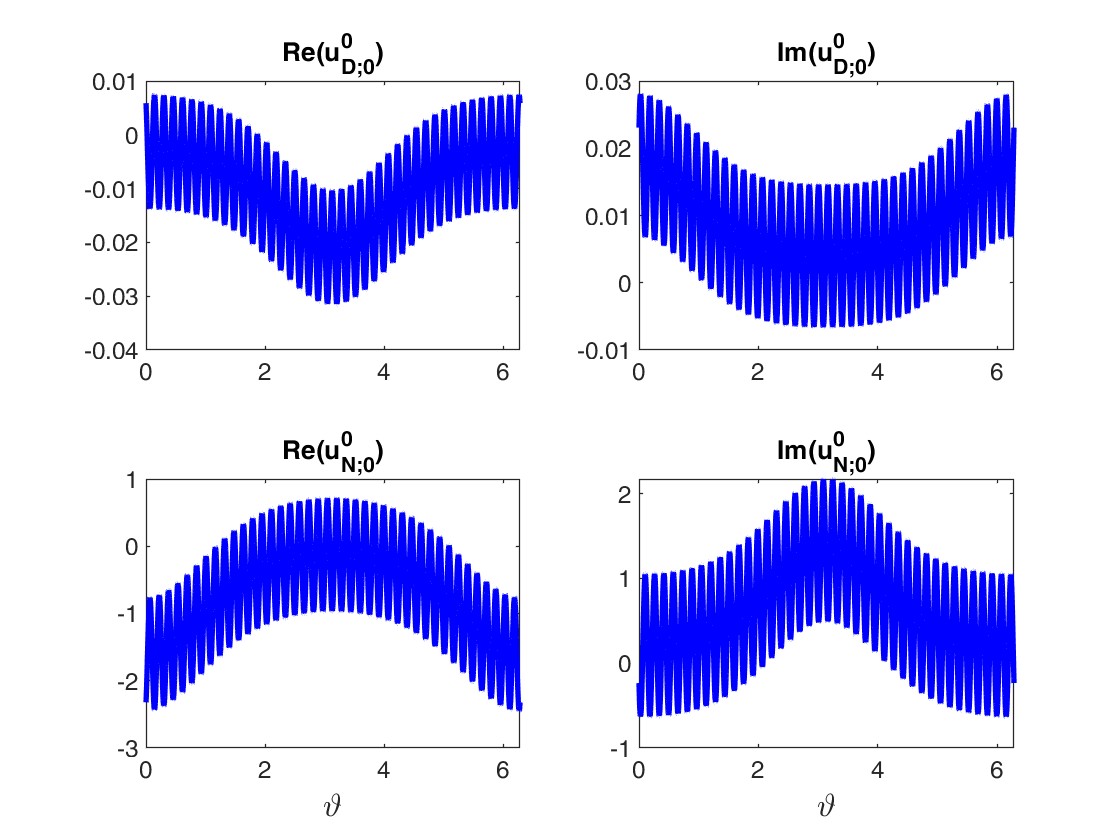}	
			
			\includegraphics[height=3cm,width=8cm]{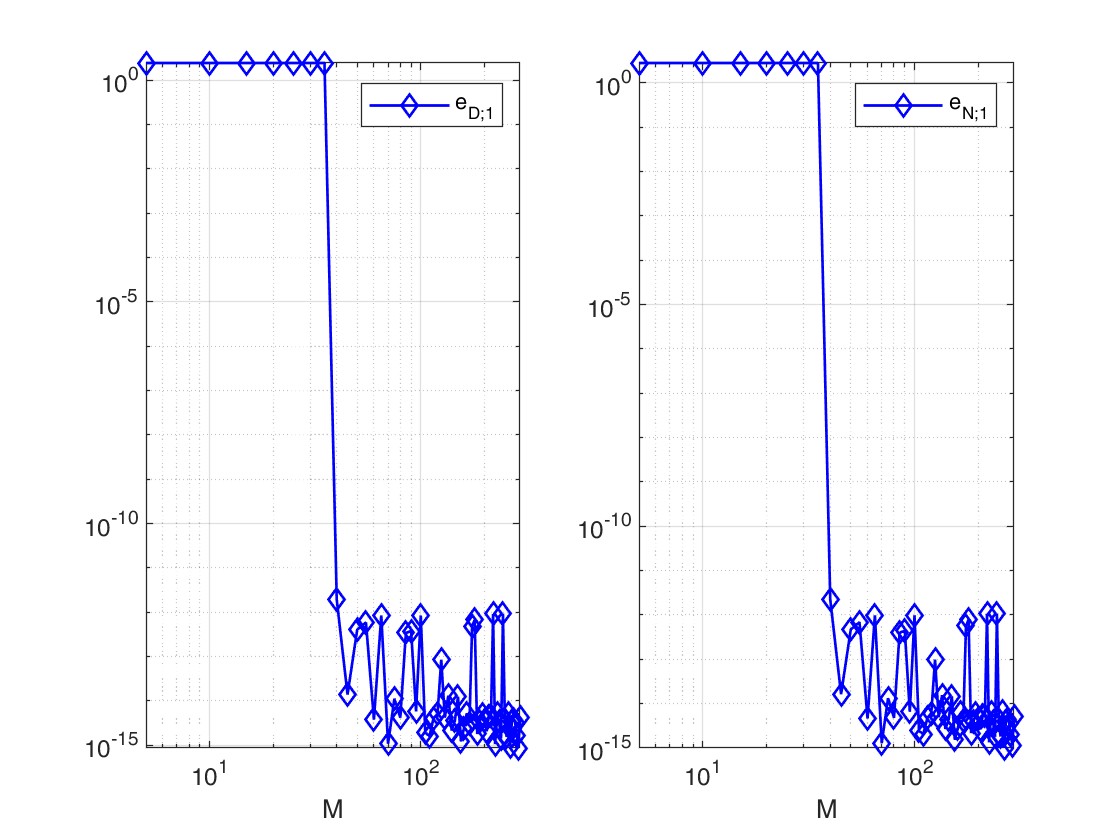}\includegraphics[height=3cm,width=8cm]{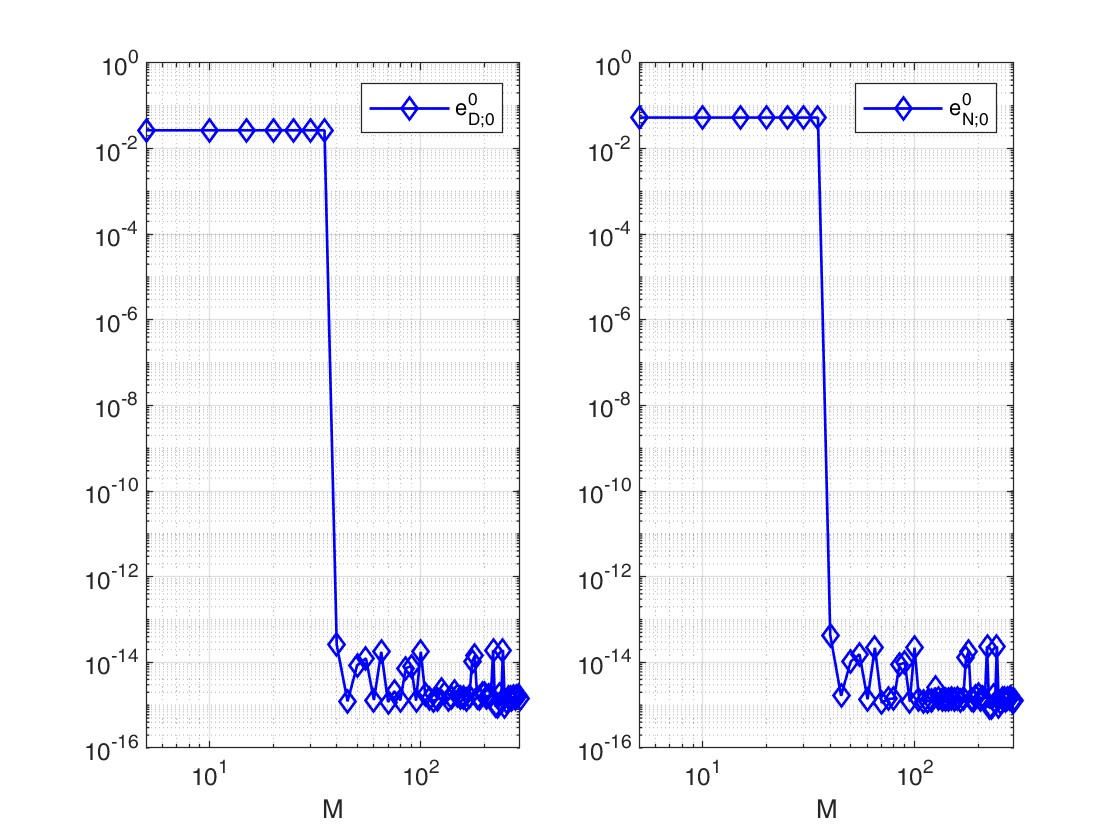}			
				\caption{Example 4. Behaviour of the traces $\uDuno, \uNuno$ (first two rows and columns) and $\uDo, \uNo$ (first two rows and second and last rows) with respect to the angle $\vartheta$. Log-log scale of the trace errors (last row).}
		\label{fig:ex4_traces}
	\end{center}
\end{figure}

In Figure \ref{fig:ex4_traces} we represent the real and imaginary parts of the four reference traces behaviour with respect to $\vartheta\in[0,2\pi]$. Also in this case the error estimates provided by Theorem \ref{th:estimates} are confirmed, as the plot of the trace errors show. Indeed, because of the regularity of the datum, we observe an exponential decay for all the errors.

\clearpage
\appendix

\section{Some estimates for Bessel, Hankel, and exponential functions}\label{sec:appendix}

In the first lemma we prove some monotonicity estimates for Bessel functions
and some bounds on their ratio.

\begin{lemma}\label{lm:lemma_A1}
\label{Propratios}For all $0<x\leq y\leq m$ and $m\in\mathbb{N}_{\geq1}$ it
holds%
\begin{align}
J_{m}\left(  x\right)   &  \leq\operatorname*{e}\nolimits^{y-x}\left(
\frac{x}{y}\right)  ^{m}J_{m}\left(  y\right)  ,\label{Jmdiffarg}\\
\left\vert \frac{\left(  H_{m}^{\left(  1\right)  }\right)  ^{\prime}\left(
	x\right)  }{H_{m}\left(  x\right)  }\right\vert  &  \leq\frac{10}{3}\left\vert
\frac{J_{m}^{\prime}\left(  x\right)  }{J_{m}\left(  x\right)  }\right\vert
\leq\frac{10}{3}\frac{m}{x}.\label{Hankelmono}%
\end{align}

\end{lemma}%
\proof
The proof of (\ref{Jmdiffarg}) follows by combining (1.9), (1.10) in
\cite{ratio_bessel_ismail}. 

Next we prove (\ref{Hankelmono}). To estimate the ratios in (\ref{Hankelmono})
we employ the general cylinder function (see, e.g., \cite{Laforgia_bessel}),
given for $\nu>0$ by:%
\[
\mathcal{C}_{\nu}\left(  x,\alpha\right)  =J_{\nu}\left(  x\right)  \cos
\alpha-Y_{\nu}\left(  x\right)  \sin\alpha,\quad0\leq\alpha<\pi.
\]
Note that%
\begin{equation}%
\begin{array}
[c]{ll}%
C_{\nu}\left(  x,0\right)  =J_{\nu}\left(  x\right)  , & C_{\nu}\left(
x,\frac{\pi}{2}\right)  =-Y_{\nu}\left(  x\right)  ,\\
C_{\nu}\left(  x,\frac{\pi}{4}\right)  =\sqrt{2}\left(  J_{\nu}\left(
x\right)  -Y_{\nu}\left(  x\right)  \right)  . & C_{\nu}\left(  x,\frac{3\pi
}{4}\right)  =-\sqrt{2}\left(  J_{\nu}\left(  x\right)  +Y_{\nu}\left(
x\right)  \right)  .
\end{array}
\label{genbessel}%
\end{equation}
Let $j_{m,1}$, $j_{m,1}^{\prime}$, $y_{m,1}$, $y_{m,1}^{\prime}$ denote the
first positive zeroes of $J_{m}$, $J_{m}^{\prime}$, $Y_{m}$, $Y_{m}^{\prime}$.
From \cite[10.2.2, 10.7.1, 10.7.4]{NIST:DLMF} it follows that%
\[%
\begin{array}
[c]{cc}%
J_{m}\left(  x\right)  \geq0\quad\forall x\in\left[  0,j_{m,1}\right]  , &
J_{m}^{\prime}\left(  x\right)  \geq0\quad\forall x\in\left[  0,j_{m,1}%
^{\prime}\right]  ,\\
Y_{m}\left(  x\right)  \leq0\quad\forall x\in\left[  0,y_{m,1}\right]  , &
Y_{m}^{\prime}\left(  x\right)  \geq0\quad\forall x\in\left[  0,y_{m,1}%
^{\prime}\right]  .
\end{array}
\]
We use \cite[10.21.3]{NIST:DLMF} for: $m\leq j_{m,1}^{\prime}<y_{m,1}%
<y_{m,1}^{\prime}<j_{m,1}$. Hence, \cite[(2.3)]{Laforgia_bessel} gives us%
\begin{equation}
\left\vert \frac{J_{m}^{\prime}\left(  x\right)  }{J_{m}\left(  x\right)
}\right\vert =\frac{J_{m}^{\prime}\left(  x\right)  }{J_{m}\left(  x\right)
}=\frac{C_{m}^{\prime}\left(  x,0\right)  }{C_{m}^{\prime}\left(  x,0\right)
}\leq\frac{m}{x}-\frac{x}{2\left(  m+1\right)  }\leq\frac{m}{x}\quad\forall
x\in\left[  0,m\right]  .\label{Jmest2}%
\end{equation}
For the ratio of the Hankel function we get%
\begin{equation}
\left\vert \frac{\left(  H_{m}^{\left(  1\right)  }\right)  ^{\prime}\left(
x\right)  }{H_{m}^{\left(  1\right)  }\left(  x\right)  }\right\vert \leq
\sqrt{\frac{\left(  J_{m}^{\prime}\left(  x\right)  \right)  ^{2}+\left(
	Y_{m}^{\prime}\left(  x\right)  \right)  ^{2}}{J_{m}^{2}\left(  x\right)
	+Y_{m}^{2}\left(  x\right)  }}\leq\frac{J_{m}^{\prime}\left(  x\right)
+Y_{m}^{\prime}\left(  x\right)  }{-Y_{m}\left(  x\right)  }\quad\forall
x\in\left[  0,m\right]  .\label{Hm1}%
\end{equation}
We use $C_{m}^{\prime}\left(  x,\frac{\pi}{4}\right)  =\sqrt{2}\left(
J_{m}^{\prime}\left(  x\right)  -Y_{m}^{\prime}\left(  x\right)  \right)
\ $(cf. \ref{genbessel}). From the asymptotic behaviour of $J_{m}^{\prime}$
and $Y_{m}^{\prime}$ as $x\rightarrow0$ (see \cite[10.2.2, 10.7.1,
10.7.4]{NIST:DLMF}) we obtain $\lim_{x\rightarrow0}C_{m}^{\prime}\left(
x,\frac{\pi}{4}\right)  =-\infty$ and%
\[
\sqrt{2}\left(  J_{m}^{\prime}\left(  x\right)  -Y_{m}^{\prime}\left(
x\right)  \right)  \leq0\quad\forall x\in\left[  0,c_{m,1}^{\prime}\right]  ,
\]
where $c_{m,1}^{\prime}$ is the first positive zero of $C_{m}^{\prime}\left(
\cdot,\frac{\pi}{4}\right)  $. From Corollary 2 in
\cite{Muldoon_zeroes_bessel_ratio}, it follows $c_{m,1}\geq m$ so that
$C_{m}^{\prime}\left(  x,\frac{\pi}{4}\right)  \leq0$ for $x\in\left[
0,m\right]  $. Hence,%
\[
J_{m}^{\prime}\left(  x\right)  \leq Y_{m}^{\prime}\left(  x\right)
\quad\forall x\in\left[  0,m\right]
\]
and in view of (\ref{Hm1}) we obtain%
\begin{equation}
\left\vert \frac{\left(  H_{m}^{\left(  1\right)  }\right)  ^{\prime}\left(
x\right)  }{H_{m}^{\left(  1\right)  }\left(  x\right)  }\right\vert
\leq2\frac{Y_{m}^{\prime}\left(  x\right)  }{-Y_{m}\left(  x\right)
}\label{Hmprimeest}%
\end{equation}
From \cite[(A.1)]{Capdebosq_ball} the first inequality in%
\[
\frac{Y_{m}^{\prime}\left(  x\right)  }{-Y_{m}\left(  x\right)  }\leq\frac
{5}{3}\left\vert \frac{J_{n}^{\prime}\left(  x\right)  }{J_{n}\left(
x\right)  }\right\vert \overset{\text{(\ref{Jmest2})}}{\leq}\frac{5}{3}%
\frac{m}{x}\quad\forall x\in\left[  0,m\right]
\]
follows. The combination with (\ref{Hmprimeest}) finally leads to%
\[
\left\vert \frac{\left(  H_{m}^{\left(  1\right)  }\right)  ^{\prime}\left(
x\right)  }{H_{m}^{\left(  1\right)  }\left(  x\right)  }\right\vert \leq
\frac{10}{3}\left\vert \frac{J_{n}^{\prime}\left(  x\right)  }{J_{n}\left(
x\right)  }\right\vert \leq\frac{10}{3}\frac{m}{x}\quad\forall x\in\left[
0,m\right]  .
\]%
\endproof

We end up this section with two technical estimates for the exponential function.

\begin{lemma}
For $0<r_{1}<r_{0}$, $\kappa_{0}>0$ and $m\in\mathbb{N}_{0}$ it holds%
\begin{equation}
\operatorname*{e}\nolimits^{\kappa_{0}\left(  r_{0}-r_{1}\right)  }\left(
\frac{r_{1}}{r_{0}}\right)  ^{m}\leq\operatorname*{e}\nolimits^{\left(
	\kappa_{0}r_{0}-m\right)  \left(  \frac{r_{0}-r_{1}}{r_{0}}\right)
}.\label{expest}%
\end{equation}
For $\eta\in\left]  0,1\right[  $, let $\mu>\operatorname*{e}/\left(
2\eta\right)  $ and set $\nu=\frac{1}{2}\left(  \frac{r_{0}-r_{1}}{r_{0}%
}\right)  \left(  \frac{\mu-1}{\mu}\right)  >0$. Then, for $m\geq\mu
\max\left\{  \kappa_{0},\kappa_{1}\right\}  r_{0}$ it holds%
\begin{equation}
\frac{\left\vert \kappa_{0}-\kappa_{1}\right\vert r_{1}}{m^{1/3}%
}\operatorname*{e}\nolimits^{\kappa_{0}\left(  r_{0}-r_{1}\right)  }\left(
\frac{r_{1}}{r_{0}}\right)  ^{m}\leq\nu^{-1}\operatorname*{e}\nolimits^{-\nu
	m}.\label{combineexpest}%
\end{equation}

\end{lemma}%
\proof
For $0<r_{1}<r_{0}$, it holds%
\begin{equation}
\operatorname*{e}\nolimits^{\kappa_{0}\left(  r_{0}-r_{1}\right)  }\left(
\frac{r_{1}}{r_{0}}\right)  ^{m}=\operatorname*{e}\nolimits^{\kappa_{0}\left(
r_{0}-r_{1}\right)  -m\left(  \log r_{0}-\log r_{1}\right)  }.\label{exp1}%
\end{equation}
A Taylor expansion around $r_{1}=r_{0}$ yields%

\[
\log r_{0}-\log r_{1}=\sum_{n=1}^{\infty}\frac{1}{n}\left(  \frac{r_{0}-r_{1}%
}{r_{0}}\right)  ^{n}\geq\frac{r_{0}-r_{1}}{r_{0}}.
\]
The combination with (\ref{exp1}) results in%
\[
\operatorname*{e}\nolimits^{\kappa_{0}\left(  r_{0}-r_{1}\right)  }\left(
\frac{r_{1}}{r_{0}}\right)  ^{m}\leq\operatorname*{e}\nolimits^{\left(
\kappa_{0}r_{0}-m\right)  \left(  \frac{r_{0}-r_{1}}{r_{0}}\right)  }.
\]

Let $\mu>1$ be as stated in the lemma. Then%
\begin{equation}
\operatorname*{e}\nolimits^{\kappa_{0}\left(  r_{0}-r_{1}\right)  }\left(
\frac{r_{1}}{r_{0}}\right)  ^{m}\leq\operatorname*{e}\nolimits^{\left(
\kappa_{0}r_{0}-m\right)  \left(  \frac{r_{0}-r_{1}}{r_{0}}\right)  }%
\leq\operatorname*{e}\nolimits^{-\left(  \frac{\mu-1}{\mu}\right)  \left(
\frac{r_{0}-r_{1}}{r_{0}}\right)  m}\leq\operatorname*{e}\nolimits^{-2\nu
m}\label{exp2nue}%
\end{equation}
for $\nu=\frac{1}{2}\left(  \frac{\mu-1}{\mu}\right)  \left(  \frac
{r_{0}-r_{1}}{r_{0}}\right)  >0$ and%
\[
\frac{\left\vert \kappa_{0}-\kappa_{1}\right\vert r_{1}}{m^{1/3}%
}\operatorname*{e}\nolimits^{\kappa_{0}\left(  r_{0}-r_{1}\right)  }\left(
\frac{r_{1}}{r_{0}}\right)  ^{m}\leq\frac{\left\vert \kappa_{0}-\kappa
_{1}\right\vert r_{1}}{m^{1/3}}\operatorname*{e}\nolimits^{-\nu m}%
\operatorname*{e}\nolimits^{-\nu m}.
\]
Elementary properties of the exponential function lead to%
\[
\frac{\left\vert \kappa_{0}-\kappa_{1}\right\vert r_{1}}{m^{1/3}}\leq
m^{2/3}\leq m\leq\nu^{-1}\left(  \nu m\right)  \leq\nu^{-1}\operatorname*{e}%
\nolimits^{\nu m}\quad\forall \, m\in{\color{black}\mathbb{N}_0}%
\]
so that finally%
\[
\frac{\left\vert \kappa_{0}-\kappa_{1}\right\vert r_{1}}{m^{1/3}%
}\operatorname*{e}\nolimits^{\kappa_{0}\left(  r_{0}-r_{1}\right)  }\left(
\frac{r_{1}}{r_{0}}\right)  ^{m}\leq\nu^{-1}\operatorname*{e}\nolimits^{-\nu
m}.
\]%
\endproof


\begin{thebibliography}{10}
	
\bibitem{Adams}
R.~Adams.
\newblock {\em {S}obolev {S}paces}.
\newblock Academic Press, N.Y., 1975.

\bibitem{Amini}
S.~Amini and S.~M. Kirkup.
\newblock Solution of {H}elmholtz equation in the exterior domain by elementary
boundary integral methods.
\newblock {\em J. Comput. Phys.}, 118(2):208--221, 1995.

\bibitem{BabuskaSauter}
I.~Babu{\v s}ka and S.~A. Sauter.
\newblock Is the pollution effect of the {FEM} avoidable for the {H}elmholtz
equation considering high wave numbers.
\newblock {\em SIAM, J. Numer. Anal.}, 34(6):2392--2423, 1997.

\bibitem{balac2021asymptotics}
S.~Balac, M.~Dauge, and Z.~Moitier.
\newblock Asymptotics for 2d whispering gallery modes in optical micro-disks
with radially varying index.
\newblock {\em IMA Journal of Applied Mathematics}, 86(6):1212--1265, 2021.

\bibitem{Turan_bessel_baricz}
A.~Baricz, S.~Ponnusamy, and S.~Singh.
\newblock Tur\'{a}n type inequalities for general {B}essel functions.
\newblock {\em Math. Inequal. Appl.}, 19(2):709--719, 2016.

\bibitem{Bouchra_Sa_Whis_I}
B.~Bensiali and S.~Sauter.
\newblock Whispering gallery modes for spherical symmetric heterogeneous
{H}elmholtz problems with piecewise smooth refractive index.
\newblock Technical Report in preparation, Preprint, 2025.

\bibitem{BoydDunster}
W.~G.~C. Boyd and T.~M. Dunster.
\newblock Uniform asymptotic solutions of a class of second-order linear
differential equations having a turning point and a regular singularity, with
an application to {L}egendre functions.
\newblock {\em SIAM J. Math. Anal.}, 17(2):422--450, 1986.

\bibitem{Capdebosq_ball}
Y.~Capdeboscq.
\newblock On the scattered field generated by a ball inhomogeneity of constant
index.
\newblock {\em Asymptot. Anal.}, 77(3-4):197--246, 2012.

\bibitem{Capdeboscq_3D}
Y.~Capdeboscq, G.~Leadbetter, and A.~Parker.
\newblock On the scattered field generated by a ball inhomogeneity of constant
index in dimension three.
\newblock In {\em Multi-scale and high-contrast {PDE}: from modelling, to
	mathematical analysis, to inversion}, volume 577 of {\em Contemp. Math.},
pages 61--80. Amer. Math. Soc., Providence, RI, 2012.

\bibitem{MonkChandlerWilde}
S.~Chandler-Wilde and P.~Monk.
\newblock Wave-{N}umber-{E}plicit {B}ounds in {T}ime-{H}armonic {S}cattering.
\newblock {\em SIAM J. Math. Anal.}, 39:1428--1455, 2008.

\bibitem{chandler2007galerkin}
S.~N. Chandler-Wilde and S.~Langdon.
\newblock A {G}alerkin boundary element method for high frequency scattering by
convex polygons.
\newblock {\em SIAM Journal on Numerical Analysis}, 45(2):610--640, 2007.

\bibitem{ChenZhou}
G.~Chen and J.~Zhou.
\newblock {\em Boundary {E}lement {M}ethods}.
\newblock Academic Press, New York, 1992.

\bibitem{Hiptmair_multiple_trace}
X.~Claeys and R.~Hiptmair.
\newblock Integral equations for acoustic scattering by partially impenetrable
composite objects.
\newblock {\em Integral Equations Operator Theory}, 81(2):151--189, 2015.

\bibitem{NIST:DLMF}
{\it NIST Digital Library of Mathematical Functions}.
\newblock http://dlmf.nist.gov/, Release 1.0.13 of 2016-09-16.
\newblock F.~W.~J. Olver, A.~B. {Olde Daalhuis}, D.~W. Lozier, B.~I. Schneider,
R.~F. Boisvert, C.~W. Clark, B.~R. Miller and B.~V. Saunders, eds.

\bibitem{ecevit2022spectral}
F.~Ecevit, Y.~Boubendir, A.~Anand, and S.~Lazergui.
\newblock Spectral {G}alerkin boundary element methods for high-frequency
sound-hard scattering problems.
\newblock {\em Numerische Mathematik}, 150(3):803--847, 2022.

\bibitem{florian2023skeleton}
F.~Florian, R.~Hiptmair, and S.~A. Sauter.
\newblock Skeleton integral equations for acoustic transmission problems with
varying coefficients.
\newblock {\em SIAM Journal on Mathematical Analysis}, 56(5):6232--6267, 2024.

\bibitem{GrHiSa_pw_Lip}
B.~Gr{\"{a}}{\textup{\ss}}le, R.~Hiptmair, and S.~A. Sauter.
\newblock Skeleton {I}ntegral {E}quations for {H}elmholtz {T}ransmission
{P}roblems: Part {II}: {P}iecewise {L}ipschitz {C}oefficients and {P}urely
{I}maginary {F}requencies.
\newblock Technical Report in preparation, Preprint, 2025.

\bibitem{Poignard_vogelius}
D.~J. Hansen, C.~Poignard, and M.~S. Vogelius.
\newblock Asymptotically precise norm estimates of scattering from a small
circular inhomogeneity.
\newblock {\em Appl. Anal.}, 86(4):433--458, 2007.

\bibitem{hiptmair2003coercive}
R.~Hiptmair.
\newblock Coercive combined field integral equations.
\newblock {\em Journal of Numerical Mathematics}, 11(2):115--134, 2003.

\bibitem{hiptmair-moiola-perugia09b}
R.~Hiptmair, A.~Moiola, and I.~Perugia.
\newblock Plane wave discontinuous {G}alerkin methods for the 2{D} {H}elmholtz
equation: analysis of the {$p$}-version.
\newblock {\em SIAM J. Numer. Anal.}, 49(1):264--284, 2011.

\bibitem{Ihlenburg}
F.~Ihlenburg and I.~Babu{\v s}ka.
\newblock Finite {E}lement {S}olution to the {H}elmholtz {E}quation with {H}igh
{W}ave {N}umber. {P}art {I}: {T}he h-version of the {FEM}.
\newblock {\em Comp. Math. Appl.}, 39(9):9--37, 1995.

\bibitem{ilchenko2006optical}
V.~S. Ilchenko and A.~B. Matsko.
\newblock Optical resonators with whispering-gallery modes-part ii:
applications.
\newblock {\em IEEE Journal of selected topics in quantum electronics},
12(1):15--32, 2006.

\bibitem{ratio_bessel_ismail}
M.~E.~H. Ismail and M.~E. Muldoon.
\newblock Certain monotonicity properties of {B}essel functions.
\newblock {\em J. Math. Anal. Appl.}, 118(1):145--150, 1986.

\bibitem{jerez2020high}
C.~Jerez-Hanckes and J.~Pinto.
\newblock High-order {G}alerkin method for {H}elmholtz and {L}aplace problems
on multiple open arcs.
\newblock {\em ESAIM: Mathematical Modelling and Numerical Analysis},
54(6):1975--2009, 2020.

\bibitem{jerez2022spectral}
C.~Jerez-Hanckes and J.~Pinto.
\newblock Spectral {G}alerkin method for solving {H}elmholtz boundary integral
equations on smooth screens.
\newblock {\em IMA Journal of Numerical Analysis}, 42(4):3571--3608, 2022.

\bibitem{Laforgia_bessel}
A.~Laforgia and P.~Natalini.
\newblock Some inequalities for modified {B}essel functions.
\newblock {\em J. Inequal. Appl.}, pages Art. ID 253035, 10, 2010.

\bibitem{MelenkLoehndorf}
M.~L\"ohndorf and J.~M. Melenk.
\newblock Wavenumber-explicit {$hp$}-{BEM} for high frequency scattering.
\newblock {\em SIAM J. Numer. Anal.}, 49(6):2340--2363, 2011.

\bibitem{MelenkSauterMathComp}
J.~M. Melenk and S.~A. Sauter.
\newblock Convergence {A}nalysis for {F}inite {E}lement {D}iscretizations of
the {H}elmholtz equation with {D}irichlet-to-{N}eumann boundary condition.
\newblock {\em Math. Comp}, 79:1871--1914, 2010.

\bibitem{MoiolaSpence_resonance_2017}
A.~{Moiola} and E.~A. {Spence}.
\newblock {Acoustic transmission problems: wavenumber-explicit bounds and
	resonance-free regions}.
\newblock {\em ArXiv e-prints}, Feb. 2017.

\bibitem{Muldoon_zeroes_bessel_ratio}
M.~E. Muldoon and R.~Spigler.
\newblock Some remarks on zeros of cylinder functions.
\newblock {\em SIAM J. Math. Anal.}, 15(6):1231--1233, 1984.

\bibitem{Nedelec01}
J.~C. N{\'e}d{\'e}lec.
\newblock {\em Acoustic and {E}lectromagnetic {E}quations}.
\newblock Springer, New York, 2001.

\bibitem{Torres_Sauter_2}
S.~Sauter and C.~Torres.
\newblock The heterogeneous {H}elmholtz problem with spherical symmetry:
{G}reen{'}s operator and stability estimates.
\newblock {\em Asymptotic Analysis}, 125(3-4):289--325, 2021.

\bibitem{Sauter2005}
S.~A. Sauter.
\newblock A {R}efined {F}inite {E}lement {C}onvergence {T}heory for {H}ighly
{I}ndefinite {H}elmholtz {P}roblems.
\newblock {\em Computing}, 78(2):101--115, 2006.

\bibitem{SaSchw1}
S.~A. Sauter and C.~Schwab.
\newblock Quadrature for hp-{G}alerkin {BEM} in {R}3.
\newblock {\em Numer. Math.}, 78(2):211--258, 1997.

\bibitem{SauterSchwab2010}
S.~A. Sauter and C.~Schwab.
\newblock {\em Boundary Element Methods}.
\newblock Springer, Heidelberg, 2010.

\bibitem{Schatz74}
A.~Schatz.
\newblock An obeservation concerning {R}itz-{G}alerkin methods with indefinite
bilinear forms.
\newblock {\em Math. Comp.}, 28:959--962, 1974.

\bibitem{vonPetersdorff89}
T.~von Petersdorff.
\newblock Boundary {I}ntegral {E}quations for mixed {D}irichlet-, {N}eumann and
{T}ransmission problems.
\newblock {\em Math. Meth. Appl. Sci.}, 11:185--213, 1989.

\bibitem{zhu-wu12b}
L.~Zhu and H.~Wu.
\newblock Preasymptotic error analysis of {CIP}-{FEM} and {FEM} for {H}elmholtz
equation with high wave number. {P}art {II}: {$hp$} version.
\newblock {\em SIAM J. Numer. Anal.}, 51(3):1828--1852, 2013.


	
\end{thebibliography}
\end{document}